\documentclass{amsart}
\usepackage{righttag}
\usepackage{epsf}
\def\hepsffile{\leavevmode\epsffile}

\theoremstyle{plain}
\newtheorem{thm}{Theorem}[subsection]
\newtheorem{thm1}{Theorem}

\newtheorem{prop}[thm]{Proposition}

\theoremstyle{definition}
\newtheorem{defin}[thm]{Definition}

\newtheorem{emf}[thm]{}
\newtheorem{rem}[thm]{Remark}

\def\id{\protect\operatorname{id}}

\def\sign{\protect\operatorname{sign}}

\def\spin{\protect\operatorname{spin}}

\def\pr{\protect\operatorname{pr}}

\def\C{{\mathbb C}}
\def\Z{{\mathbb Z}}
\def\Q{{\mathbb Q}}
\def\R{{\mathbb R}}
\def\N{{\mathbb N}}

\def\1{\hbox{\rm\rlap {1}\hskip.03in{\rom I}}}
\def\Bbbone{{\rm1\mathchoice{\kern-0.25em}{\kern-0.25em}
	{\kern-0.2em}{\kern-0.2em}I}}

\def\p{\partial}

\begin{document}
\hyphenation{Ca-m-po}
\title[Vassiliev Invariants of Legendrian, Transverse, and Framed Knots]
{Vassiliev Invariants of Legendrian, Transverse, and Framed Knots in Contact $3$-manifolds}
\author[V.~Tchernov]{Vladimir Tchernov}
\address{Max-Planck-Institut f\"ur Mathematik, P.~O.~Box 7280, D-53072 Bonn,
Germany}
\email{chernov@mpim-bonn.mpg.de}
\begin{abstract}
We show that for a large class of contact $3$-manifolds 
the groups of Vassiliev
invariants of Legendrian and 
of framed knots are canonically isomorphic.
As a corollary, we obtain that the group of finite order 
Arnold's $J^+$-type invariants of wave fronts on a surface $F$ is isomorphic
to the group of Vassiliev invariants of framed knots in the spherical
cotangent bundle $ST^*F$ of $F$.

On the other hand we construct the first examples of contact manifolds 
for which Vassiliev invariants of Legendrian knots can
distinguish Legendrian knots that realize isotopic framed knots and are
homotopic as Legendrian immersions.
\end{abstract}
\maketitle

\section{Introduction} 
In this section we describe the main results of the paper. 
(In case any of the terminology appears to be new to the reader, the
corresponding definitions are given in the next section.)

If a contact structure on a $3$-manifold is cooriented, then every
Legendrian knot (i.e. a knot that is everywhere tangent to the contact
distribution) has a natural framing (a continuous normal vector field).
Hence when studying Legendrian knots in such contact manifolds the main
question is to distinguish those of them that realize isotopic framed knots.

Similarly if the contact structure is parallelized, then every transverse
knot (i.e. a knot that is everywhere transverse to the contact distribution)
also has a natural framing, and when studying transverse knots in such
contact manifolds again the main question is to distinguish those of them that 
realize isotopic framed knots.

Vassiliev invariants proved to be an extremely useful tool in the study of
framed knots, and the conjecture is that they are sufficient to distinguish all
the isotopy classes of framed knots. Vassiliev invariants can also be easily defined in the categories of
Legendrian and of transverse knots.
In this paper we study the relationship between the groups of Vassiliev
invariants of these three categories of knots, and explore when these 
invariants can be used to distinguish Legendrian knots that
realize isotopic framed knots. 

Consider a contact manifold $M$ with a cooriented contact structure.
Fix an Abelian group $\mathcal A$, a connected component $\mathcal F$ 
of the space of framed immersions of $S^1$ into $M$, and a connected component
$\mathcal L\subset \mathcal F$ of the space of
Legendrian immersions of $S^1$ into $M$.
We study the relation between the groups of $\mathcal A$-valued 
Vassiliev invariants of
framed knots from $\mathcal F$ and of $\mathcal A$-valued 
Vassiliev invariants of Legendrian knots from $\mathcal L$. The main results
obtained in this paper are described below.

\begin{thm1}\label{first}
The groups of $\mathcal A$-valued 
Vassiliev invariants of Legendrian knots from $\mathcal L$ and of framed knots 
from $\mathcal F$ are canonically
isomorphic,
provided that the Euler class of the
contact bundle vanishes on every $\alpha\in H_2(M,\Z)$ realizable by a
mapping $\mu:S^1\times S^1\rightarrow M$.
\end{thm1}

(See Theorem~\ref{isomorphismobtained} and
Proposition~\ref{interpretationconditionII}.)

Using Theorem~\ref{first} we show that:

\begin{thm1}
The groups of $\mathcal A$-valued
Vassiliev invariants of Legendrian knots from $\mathcal L$ and of framed
knots
from $\mathcal F$ are canonically
isomorphic, provided that one of the following conditions holds:
\begin{description}
\item[1] the contact structure is tight;
\item[2] the Euler class of the contact bundle is in the torsion of
$H^2(M,\Z)$ (in particular if the Euler class is zero).
\item[3] the contact manifold is closed and admits a metric of negative
sectional curvature.
\end{description}
\end{thm1}

(See Sections~\ref{homologysphere}
and~\ref{negativecurvature} and Theorem~\ref{tight}.) 

As a corollary, we get that for any surface $F$ 
the group of finite order Arnold's $J^+$-type 
invariants of wave fronts on $F$ is isomorphic to the group of Vassiliev invariants 
of framed knots in the spherical cotangent bundle $ST^*F$ of $F$.

Previously the isomorphism of the groups of Vassiliev invariants of
Legendrian and of framed knots was known only in the 
case where $\mathcal A=\C$ and $M$ is the standard contact $\R^3$
(result of D.~Fuchs and S.~Tabachnikov~\cite{FuchsTabachnikov})
or the standard contact solid-torus (result of J.~Hill~\cite{Hill}).
The proofs of these isomorphisms were based on the fact that 
for the $\C$-valued Vassiliev invariants of framed knots in these manifolds 
there exists a universal Vassiliev invariant also known as the Kontsevich
integral. (Currently the existence of the Kontsevich integral is known only for 
a total space of an $\R^1$-bundle over a compact oriented surface
with boundary, see the paper~\cite{AMR} of Andersen, Mattes, and
Reshetikhin.)

Thus the approach used in~\cite{FuchsTabachnikov} and~\cite{Hill} 
to show the isomorphism of
the groups of Vassiliev invariants 
is not applicable for almost all contact $3$-manifolds
and Abelian groups $\mathcal A$ and our results appear to be a strong
generalization of the results of Fuchs, Tabachnikov and Hill.

We also construct the first examples where Vassiliev invariants can be used
to distinguish Legendrian knots that realize isotopic framed knots and are
homotopic as Legendrian immersions. These are also the first examples 
where the groups of
Vassiliev invariants of Legendrian and of framed knots from the
corresponding components of the spaces of Legendrian and of framed immersions
are not canonically isomorphic.

\begin{thm1}
The manifold  $S^1\times S^2$ admits infinitely many cooriented contact structures
for which there exist Legendrian knots that can be distinguished by
$\Z$-valued 
Vassiliev invariants even though they realize isotopic framed knots and are
homotopic as Legendrian immersions.
\end{thm1}

(See Theorem~\ref{example1} and Theorem~\ref{example2} in which the similar
result is proved for any orientable total space of an $S^1$-bundle over a nonorientable
surface of a sufficiently high genus.)

For transverse knots we obtain the following result
(see Theorem~\ref{isomorphismobtainedtransverse}):

\begin{thm1}
Let $(M,C)$ be a contact manifold with a parallelized contact structure,
then the groups of $\mathcal A$-valued Vassiliev invariants of
transverse and of framed knots (from the corresponding components of the
spaces of transverse and of framed immersions) are canonically isomorphic. 
\end{thm1}

\section{Conventions and definitions.}\label{definitions}
In this paper $\mathcal A$ is an Abelian group
(not necessarily torsion free), and $M$ is a 
connected oriented $3$-dimensional Riemannian manifold (not necessarily
compact).

A {\em contact structure\/} on a $3$-dimensional manifold $M$ is a smooth
field $\{ C_x\subset T_xM|x\in M\}$ of tangent $2$-dimensional planes,
locally defined as a kernel of a differential $1$-form $\alpha$ with non-vanishing
$\alpha\wedge d\alpha$. A manifold with a contact structure possesses
the canonical orientation determined by the volume form $\alpha\wedge\ d
\alpha$. The standard contact structure in $\R^3$ is the kernel of the $1$-form
$\alpha=ydx-dz$.

A {\em contact element\/} on a manifold is a hyperplane
in the tangent space to the manifold at a point.
For a surface $F$ we denote by $ST^*F$ the space of all
cooriented (transversally oriented) contact elements of $F$. This space is
the spherical cotangent bundle of $F$. Its natural contact structure is the
distribution of tangent hyperplanes given by the condition
that the velocity vector of the incidence point of a contact element
belongs to the element.

A contact structure is {\em cooriented\/} if the
$2$-dimensional planes defining the contact structure are continuously
cooriented (transversally oriented). 
A contact structure is {\em oriented\/} if the $2$-dimensional planes defining the
contact structure are continuously oriented. Since every contact manifold
has a natural orientation we see that every cooriented contact structure is
naturally oriented and every oriented contact structure is naturally
cooriented.
A contact structure is {\em
parallelizable\/} ({\em parallelized\/}) if the $2$-dimensional vector bundle 
$\{ C_x \}$ over $M$ is trivializable (trivialized). Since every contact manifold
has a canonical orientation, one can see that every parallelized contact
structure is naturally cooriented. A contact structure $C$ on a manifold $M$
is said to be {\em overtwisted\/} if there exists a $2$-disk $D$ embedded
into $M$ such that the boundary $\p D$ is tangent to $C$ while the disk $D$
is transverse to $C$ along $\p D$. Not overtwisted contact structures are
called {\em tight\/}.

A {\em curve\/} in $M$ is an immersion of $S^1$ into $M$. 
(All curves have the natural orientation induced by the orientation of
$S^1$.) A {\em framed curve\/} in $M$ is a curve in $M$ equipped 
with a continuous unit normal vector field.

A {\em Legendrian curve\/} in a contact manifold $(M,C)$ is a curve
in $M$ that is everywhere tangent to $C$. If the contact structure on
$M$ is cooriented, then every Legendrian curve has a natural framing given 
by the unit normals to the planes of the contact structure that point in the
direction specified by the coorientation.

To a Legendrian curve $K_l$ in a contact manifold with a parallelized 
contact structure one can associate an integer that is the
number of revolutions of the direction of the velocity vector of $K_l$ (with
respect to the chosen frames in $C$) under traversing $K_l$ according
to the orientation. This integer is called the {\em Maslov number\/} of
$K_l$. The set of Maslov numbers enumerates the
set of the connected components of the space of Legendrian 
curves in $\R^3$ (cf.~\ref{h-principleLegendrian}).

A {\em transverse\/} curve in a contact manifold $(M,C)$ is a curve in $M$
that is everywhere transverse to $C$. If the
contact structure on $M$ is parallelized, then a transverse curve has a
natural framing given by the unit normals corresponding to  
the projections  of the first of the two coordinate
vectors of the contact planes on the $2$-planes orthogonal 
to the velocity vectors of the curve. 
A transverse curve in a contact manifold
with a cooriented contact structures is said to be {\em positive\/} if at
every point the velocity vector of the curve points into the coorienting
half-plane, and it is said to be {\em negative\/} otherwise. There are two
connected components of the space of transverse curves in $\R^3$, they
consist of positive and negative transverse curves respectively. In general
if $(M,C)$ is a contact manifold with a cooriented contact structure, then every
connected component of the space of unframed curves contains two connected
components of the space of transverse curves. They consist of positive and
negative transverse curves respectively.

A {\em knot (framed knot)\/} in $M$ is an embedding (framed embedding) of $S^1$ into $M$.
In a similar way we define Legendrian and transverse knots in $M$.

A {\em singular (framed)\/} knot with $n$ double points is a curve (framed curve)
in $M$ whose only singularities are $n$ transverse double points.
An {\em isotopy\/} of a singular (framed) knot 
with $n$ double points is a path in the space of singular (framed) knots with
$n$ double points under which the preimages of the double points on $S^1$
change continuously.

An $\mathcal A$-valued framed (resp. Legendrian, resp. transverse) knot invariant is 
an $\mathcal A$-valued function on the set of the isotopy classes of framed (resp.
Legendrian, resp. transverse) knots.

A transverse double point $t$ of a singular knot can be resolved in two 
essentially different ways. We say that a resolution of a double point is
positive (resp. negative) if the tangent vector to the
first strand, the tangent vector to the second strand, and the vector from
the second strand to the first form the positive $3$-frame. (This does 
not depend on the order of the strands).
If the singular knot is Legendrian (resp. transverse), 
then these resolution can be made in the
category of Legendrian (resp. transverse) knots.

A singular framed (resp. Legendrian, resp. transverse) knot $K$ with $(n+1)$ 
transverse double points
admits $2^{n+1}$ possible resolutions of the double points. The sign of the resolution 
is put to be $+$ if the number of negatively resolved double points is even, and
it is put to be $-$ otherwise. 
Let $x$ be an $\mathcal A$-valued invariant of framed (resp. Legendrian,
resp. transverse) knots. The invariant $x$ is said to be of {\em finite
order\/} (or {\em Vassiliev invariant\/}) if there exists a nonnegative 
integer $n$ such that for any singular knot $K_s$ with $(n+1)$
transverse double points the sum (with appropriate signs) of the values of $x$ on the nonsingular
knots obtained by the $2^{n+1}$ resolutions of the double points is zero. 
An invariant is said to be of order not greater than $n$ (of order $\leq n$) if $n$
can be chosen as the integer in the definition above. The group of $\mathcal
A$-valued finite order invariants has an increasing filtration by the
subgroups of the invariants of order $\leq n$.

\begin{emf}\label{h-principleLegendrian}
{\em $h$-principle for Legendrian curves.\/}
For $(M,C)$ a contact manifold with a cooriented contact structure, we put
$CM$ to be the total space of the fiberwise spherization of the contact
bundle, and we put $\pr:CM\rightarrow M$ to be the corresponding locally
trivial $S^1$-fibration. The $h$-principle proved for the Legendrian curves
by M.~Gromov (\cite{Gromov}, pp.338-339) says that the space of Legendrian
curves in
$(M,C)$ is weak homotopy equivalent to the space of free loops
$\Omega CM$
in $CM$. The equivalence is given by mapping a point of a Legendrian curve
to the point of $CM$ corresponding to the direction of the velocity vector
of the curve at this point. In particular the $h$-principle implies that the
set of the connected components of the space of Legendrian curves can be
naturally identified with the set of the conjugacy classes of elements of
$\pi_1(CM)$.
\end{emf}

\begin{emf}\label{description}{\em Description of Legendrian and of transverse
knots in $\R^3$.\/} 
The contact Darboux theorem says that every contact $3$-manifold $(M,C)$ is
locally contactomorphic to 
$\R^3$ with the standard contact structure that is the kernel of the $1$-form 
$\alpha=ydx-dz$. A chart in which $(M,C)$ is contactomorphic to the standard
contact $\R^3$ is called {\em a Darboux chart.\/}

Transverse and Legendrian knots in the standard contact $\R^3$
are conveniently presented by the projections 
into the plane $(x,z)$. Identify a  point $(x, y, z)\in \R^3$ with the
point $(x,z)\in \R^2$ furnished with the fixed  
direction of an unoriented 
straight line through $(x,z)$ with the slope $y$. Then the curve in $\R^3$
is a one parameter family of points with (non-vertical) directions in $\R^2$. 

A curve in
$\R^3$ is transverse if and only if the corresponding curve in $\R^2$ is
never tangent to the chosen directions along itself.

While a generic regular curve has a regular projection into the
$(x,z)$-plane, the projection of a generic Legendrian curve into the
$(x,z)$-plane has isolated critical points (since all the planes of the
contact structure are parallel to the $y$-axis). Hence the projection of a
generic Legendrian curve may have cusps. A curve in $\R^3$ is Legendrian if
and only if the corresponding planar curve with cusps 
is everywhere tangent to the field of directions. In particular this field
is determined by the curve with cusps.
\end{emf}

\section{Isomorphisms of the groups of Vassiliev invariants of Legendrian,
of transverse, and of framed knots.}

\subsection{Isomorphism between the groups of order $\leq n$ invariants of
Legendrian and of framed knots}

Let $(M,C)$ be a contact manifold with a cooriented contact structure.
Let $\mathcal L$ be a connected component of
the space of Legendrian curves in $M$, and let $\mathcal F$ be the
connected
component of the space of framed curves that contains $\mathcal L$. 
(Such a
component exists because a Legendrian curve in a manifold with a
cooriented contact structure is naturally framed, and a path
in the space of Legendrian curves corresponds to a path in the space of
framed curves.)
Let $V_n^{\mathcal L}$ (resp. $W_n^{\mathcal F}$) 
be the group of $\mathcal A$-valued order $\leq n$ invariants of Legendrian
(resp. framed) knots from $\mathcal L$ (resp.
from $\mathcal F$). Clearly every invariant $y\in W_n^{\mathcal F}$ restricted to
the category of Legendrian knots in $\mathcal L$ is an element $\phi(y)\in
V_n^{\mathcal L}$. This gives a homomorphism $\phi:W_n^{\mathcal
F}\rightarrow V_n^{\mathcal L}$.

\begin{thm}\label{isomorphism} Let $(M, C)$ be a contact manifold with a 
cooriented contact structure. Let $\mathcal L$ be a connected component of
the space of Legendrian curves in $M$, and let $\mathcal F$ be the
connected
component of the space of framed curves that contains $\mathcal L$.
Then the following two statements {\bf a} and {\bf b} are equivalent.
\begin{description}
\item[a] $x(K_1)=x(K_2)$ 
for any $x\in V_n^{\mathcal L}$ and any knots $K_1, K_2\in
\mathcal L$ representing isotopic framed knots. 
\item[b] $\phi:W_n^{\mathcal F}\rightarrow V_n^{\mathcal L}$ is a canonical
isomorphism.
%
\end{description}
\end{thm}

If the mapping from the isotopy classes of Legendrian knots
in $\mathcal L$ to the isotopy classes of framed knots in $\mathcal F$ is
surjective, then the proof of Theorem~\ref{isomorphism} is obvious.
However in general this mapping is not surjective and the
proof of Theorem~\ref{isomorphism} is given in
Section~\ref{proofisomorphism}. 
The famous Bennequin inequality shows that this mapping is not surjective
even in the case where $M$ is the standard contact $\R^3$.


Theorem~\ref{isomorphism} implies that 
to obtain the isomorphism between the groups $W_n^{\mathcal F}$ and
$V_n^{\mathcal L}$ it suffices to show that statement $\bf a$ of
Theorem~\ref{isomorphism} is true for the connected components $\mathcal L$
and $\mathcal F$ of the spaces of Legendrian and of framed curves. 

\begin{emf}\label{conditions}{\em Condition $(*)$.}
In~\cite{FuchsTabachnikov} D.~Fuchs and S.~Tabachnikov 
showed that statement $\bf a$ holds for 
all the connected components of the space of Legendrian curves 
when the ambient manifold is the standard contact $\R^3$ 
and the group $\mathcal A$ is $\C$. (One can verify that the proof of this
Theorem of Fuchs and Tabachnikov goes through for $\mathcal A$ being any
Abelian group.)
They later observed~\cite{FuchsTabachnikovprivate} that since their proof 
of this fact is mostly local, the similar
fact should be true for a big class of contact manifolds. 

However in fact the proof of their theorem is not
completely local and is also based on the existence of well-defined
Bennequin invariant and Maslov number for a Legendrian knot in $\R^3$. 
In general the Bennequin invariant is not well-defined unless the knot is
zero-homologous and the Maslov number is not well-defined unless either the
knot is zero-homologous or the contact structure is parallelizable. Thus the
generalization of this Theorem to the case of manifolds other than $\R^3$
meets certain difficulties. (And in fact the corresponding result does not
hold for a big class of contact manifolds, see
Section~\ref{nonisomorphicgroups}.)  

By analyzing the
proof of the Theorem of Fuchs and Tabachnikov (see
Section~\ref{reasonstoappear}) we get that it can be
generalized 
to the case of an arbitrary contact $3$-manifold with a cooriented contact
structure, provided that the connected component $\mathcal F$
(containing $\mathcal L$) satisfies
the following 

\smallskip
\noindent {\bf Condition $(*)$:}
the connected component $\mathcal F$
of the space of framed curves contains infinitely many components of the
space of Legendrian curves. (See Proposition~\ref{interpretationconditionII} for
the homological interpretation of condition $(*)$.)
\end{emf}

This generalization of the result of Fuchs and Tabachnikov 
and Theorem~\ref{isomorphism} imply the following Theorem. 

\begin{thm}\label{isomorphismobtained}
Let $(M, C)$ be a contact manifold with a
cooriented contact structure, and let $\mathcal L$ be a connected component of
the space of Legendrian curves in $M$. Let $\mathcal F$ be the
connected component of the space of framed curves that contains $\mathcal L$.
Let $V_n^{\mathcal L}$ (resp. $W_n^{\mathcal F}$) 
be the group of $\mathcal A$-valued order $\leq n$ invariants of Legendrian
(resp. framed) knots from $\mathcal L$ (resp. 
from $\mathcal F$). 
Then the groups $V_n^{\mathcal L}$ and $W_n^{\mathcal F}$ are canonically
isomorphic, provided that $\mathcal F$ satisfies condition $(*)$.
\end{thm}

Now we give a homological interpretation of condition $(*)$.

\begin{prop}\label{interpretationconditionII}
Let $(M,C)$ be a 
contact manifold with a cooriented contact
structure, let $\chi_C\in H^2(M)$ be the Euler class of the contact bundle, 
and let $\mathcal F$ be a component of the space of framed curves
in $M$. Then $\mathcal F$ does not satisfy condition $(*)$ if and only if 
there exists $\alpha\in H_2(M, \Z)$ such that $\chi_C(\alpha)\neq 0$
and $\alpha$ is realizable by a mapping 
$\mu:S^1\times S^1\rightarrow M$ with the property that
$\mu\big|_{1\times S^1}$ is a loop free homotopic to loops realized 
by curves from $\mathcal F$.
\end{prop}

For the Proof of Proposition~\ref{interpretationconditionII} see
Subsection~\ref{ProofinterpretationconditionII}.

\begin{rem}\label{homologysphere}{\em Some immediate corollaries of
Theorem~\ref{isomorphismobtained} and the generalization of the Theorem of
Fuchs and Tabachnikov about the isomorphism of the groups of the $\C$-valued
Vassiliev invariants in the case of $M=\R^3$.}
Proposition~\ref{interpretationconditionII} implies that if the contact
structure is parallelizable (and hence the Euler class of the contact bundle
is zero) then all the connected components of the space of framed curves
satisfy condition $(*)$. Applying Theorem~\ref{isomorphismobtained} we
conclude that for any Abelian group $\mathcal A$ and 
for every connected component of the space of Legendrian curves
$\mathcal L$ and for the containing it component of the space of framed curves
$\mathcal F$ the groups $V_n^{\mathcal L}$ and $W_n^{\mathcal F}$ 
of $\mathcal A$-valued Vassiliev invariants are canonically isomorphic.

Clearly the value of the Euler class of the contact bundle is zero if $M$ is
an integer homology sphere. Hence for any Abelian group $\mathcal A$ we
obtain the isomorphism of the groups $V_n^{\mathcal L}$                     
and $W_n^{\mathcal F}$ of $\mathcal A$-valued 
Vassiliev invariants.
This generalizes the Theorem of D.~Fuchs and
S.~Tabachnikov~\cite{FuchsTabachnikov} saying that 
for the standard contact $\R^3$ 
and for $\mathcal A=\C$ the quotient groups
$V_n^{\mathcal L}/V_{n-1}^{\mathcal L}$ and $W_n^{\mathcal
F}/W_{n-1}^{\mathcal F}$ are canonically isomorphic. 

The proof of this
Theorem of Fuchs and Tabachnikov was based on the fact that for the
$\C$-valued Vassiliev invariants of framed knots in
$\R^3$ there exists the universal Vassiliev invariant
constructed by T.~Q.~T.~Le and J.~Murakami~\cite{LeMurakami}.
(For unframed knots in $\R^3$ the construction of the universal Vassiliev
invariant is the classical result of M.~Kontsevich~\cite{Kontsevich}, and
the invariant itself is the famous Kontsevich integral.)
The existence of the universal Vassiliev invariant is currently known only for a
very limited collection of $3$-manifolds, and only for $\mathcal A$ being
$\C$, $\R$, or $\Q$.
(Andersen, Mattes, Reshetikhin~\cite{AMR} proved its existence in the case
where 
$\mathcal A=\C$ and $M$ is the total space of an $\R^1$-bundle 
over a compact oriented surface $F$ with $\p F\neq \emptyset$.)

Thus the approach used in~\cite{FuchsTabachnikov} to show the isomorphism of
the quotient groups is not applicable for almost all contact $3$-manifolds
and Abelian groups $\mathcal A$, and
Theorem~\ref{isomorphismobtained} appears to be a strong
generalization of the result of Fuchs and Tabachnikov.

\end{rem}

\begin{rem}
Let $(M,C)$ be a contact manifold with a cooriented contact structure, and
let $\mathcal F$ be a connected component of the space of framed curves in
$M$. 
Theorem~\ref{isomorphismobtained} 
implies that the group of $\mathcal A$-valued order 
$\leq n$ invariants of Legendrian knots from a connected component $\mathcal
L\subset \mathcal F$
of the space of Legendrian curves does not depend on the
choice of a cooriented contact structure, 
provided that for this choice $\mathcal F$ satisfies condition~$(*)$. And hence in these cases the group can not be
used to distinguish cooriented contact structures on $M$. (See
Remark~\ref{homologysphere} and Theorems~\ref{atoroidal} and~\ref{tight} 
for the list of cases when
the connected components of the space of framed curves are known to 
satisfy condition~$(*)$.)
\end{rem}



%

\begin{emf}\label{ArnoldsJ+}{\bf Finite order Arnold's $J^+$-type
invariants of wave fronts on surfaces.\/}
A very interesting class of contact manifolds satisfying the conditions of
Theorem~\ref{isomorphismobtained} is formed by the spherical
cotangent bundles $ST^*F$ of surfaces $F$ with the natural contact
structure on $ST^*F$ (see~\ref{definitions}). The theory of the invariants 
of Legendrian knots in $ST^*F$ is often referred to as the theory of 
Arnold's~\cite{Arnoldsplit} $J^+$-type invariants of fronts on a surface $F$. 
The natural contact structure on $ST^*F$ 
is cooriented. (The coorientation is induced from the
coorientation of the contact elements of $F$.) 
One can verify that for orientable $F$ 
the standard contact structure on $ST^*F$ is parallelizable, and hence all
the components of the space of framed curves satisfy condition~$(*)$.
If $F$ is not orientable, then the standard cooriented contact structure on
$ST^*F$ is not parallelizable, but one can still verify 
(cf. Proposition 8.2.4~\cite{Tchernov}) that every connected component of the
space of framed curves satisfies condition~$(*)$. 
Hence for any Abelian group $\mathcal A$ and for 
any surface $F$ we obtain the canonical isomorphism of   
the groups of $\mathcal A$-valued order $\leq n$ invariants of Legendrian 
and of framed knots (from the corresponding components of the spaces of 
Legendrian and of framed curves in $ST^*F$ with the standard contact structure). 
Or equivalently we get that the groups of $\mathcal A$-valued order 
$\leq n$ $J^+$-type invariants of fronts on $F$
and of $\mathcal A$-valued order $\leq n$ invariants of framed
knots in $ST^*F$ (from the corresponding components of the two spaces) are canonically
isomorphic.

Previously it was known that for $F=\R^2$ and $\mathcal
A=\C$ the quotient groups 
$V_n^{\mathcal L}/V_{n-1}^{\mathcal L}$ and $W_n^{\mathcal
F}/W_{n-1}^{\mathcal F}$ are canonically isomorphic. The proof of this result of J.~W.~Hill~\cite{Hill}
was based
on the fact that for the $\C$-valued Vassiliev invariants of framed knots in
$ST^*\R^2$ there exists the universal Vassiliev invariant 
constructed by V.~Goryunov~\cite{Goryunov}.
(For unframed knots in $\R^3$ the existence of the universal Vassiliev 
invariant is the classical result of M.~Kontsevich~\cite{Kontsevich}, and
the invariant itself is the famous Kontsevich integral.)
Our results generalize the result of J.~W. Hill (even in the case of
$M=ST^*\R^2)$.
\end{emf}

The following Theorem describes another big class of contact manifolds for
which the groups of Vassiliev invariants of Legendrian and of framed knots
(from the corresponding components of the two spaces of curves) are canonically
isomorphic.

\begin{thm}\label{atoroidal}
Let $(M,C)$ be a contact manifold (with a cooriented contact
structure) such that $\pi_2(M)=0$, and for every mapping
$\mu:T^2\rightarrow M$ of the two-torus 
the homomorphism $\mu_*:\pi_1(T^2)\rightarrow \pi_1(M)$ is not injective.
Then all the components of the space of framed curves in $M$ satisfy
condition~$(*)$, and hence the groups of $\mathcal
A$-valued order $\leq n$ invariants of Legendrian and of framed knots (from the
corresponding components of the spaces of Legendrian and framed curves) 
are canonically isomorphic.
\end{thm} 

For the Proof of Theorem~\ref{atoroidal} see
Subsection~\ref{Proofatoridal}.

\begin{emf}\label{negativecurvature}
{\em The isomorphism of the 
groups of Vassiliev invariants in the case of closed manifolds admitting a
metric of negative sectional curvature and other 
corollaries of Theorem~\ref{atoroidal}.\/} 

Let $M$ be a closed manifold admitting a metric of negative sectional
curvature.
A well-known Theorem by A.~Preissman (see~\cite{Docarmo} pp. 258-265) 
says that every nontrivial commutative subgroup of the fundamental group of a 
closed $3$-dimensional manifold of negative sectional curvature is infinite 
cyclic. Hence for every mapping $\mu:T^2\rightarrow M$ the kernel of 
$\mu_*:\pi_1(T^2)=\Z\oplus\Z\rightarrow \pi_1(M)$ is nontrivial.
It is also known that the universal covering of such $M$ is diffeomorphic 
to $\R^3$, and hence $\pi_2(M)=0$. Thus every closed manifold $M$ 
admitting a metric of negative  sectional 
curvature satisfies all the conditions of Theorem~\ref{atoroidal} and for an arbitrary 
cooriented contact structure on such $M$ we obtain the isomorphism of the groups 
of $\mathcal A$-valued order $\leq n$ invariants of Legendrian and of framed
knots from the corresponding components of the spaces of Legendrian and of
framed curves.

%
\end{emf}

Another important class of contact manifolds for which every connected
component of the space of framed curves satisfies condition~$(*)$ is
formed by contact manifolds with a tight contact structure.

The following Theorem appeared as a result of discussions of 
Stefan Nemirovski and the author.

\begin{thm}\label{tight} 
Let $(M,C)$ be a contact manifold with a tight cooriented contact
structure. 
Then all the components of the space of framed curves in $M$ satisfy
condition~$(*)$, and hence the groups of $\mathcal
A$-valued order $\leq n$ invariants of Legendrian and of framed knots (from the
corresponding components of the spaces of Legendrian and framed curves) are canonically isomorphic.
\end{thm}

For the Proof of Theorem~\ref{tight} see Subsection~\ref{Prooftight}.

\subsection{Isomorphisms between the groups of Vassiliev invariants of
transverse and of framed knots}
Let $M$ be a contact manifold with a
parallelized contact structure $C$. Let $\mathcal T$ be a connected component of
the space of transverse curves in $(M,C)$, and let $\mathcal F$ be the connected
component of the space of framed curves that contains $\mathcal T$. (Such a
component exists because a transverse curve in a manifold with a
parallelized contact structure is naturally framed, and a path
in the space of transverse curves corresponds to a path in the space of
framed curves.) Let $V_n^{\mathcal T}$ (resp. $W_n^{\mathcal F}$) 
be the group of $\mathcal A$-valued order $\leq n$ invariants of transverse
(resp. framed) knots from $\mathcal T$ (resp. from $\mathcal F$). Clearly every 
invariant $y\in W_n^{\mathcal F}$ restricted to
the category of transverse knots in $\mathcal T$ is an element $\phi(y)\in
V_n^{\mathcal T}$. This gives a homomorphism $\phi:W_n^{\mathcal
F}\rightarrow V_n^{\mathcal T}$.

\begin{thm}\label{isomorphismtransverse} Let $(M, C)$ be a contact manifold 
with a parallelized contact structure. Let $\mathcal T$ be a connected component of
the space of transverse curves in $(M,C)$, and let $\mathcal F$ be the
component of the space of framed curves that contains $\mathcal T$.
Then the following two statements {\bf a} and {\bf b} are equivalent.
\begin{description}
\item[a] $x(K_1)=x(K_2)$ 
for any $x\in V_n^{\mathcal T}$ and any knots $K_1, K_2\in
\mathcal T$ representing isotopic framed knots. 
\item[b] 
$\phi:W_n^{\mathcal F}\rightarrow V_n^{\mathcal T}$ is a canonical isomorphism.
\end{description}
\end{thm} 

The proof of Theorem~\ref{isomorphismtransverse} is analogous to the Proof 
of Theorem~\ref{isomorphism}. 

Similar to the case of Theorem~\ref{isomorphism} the proof of
Theorem~\ref{isomorphismtransverse} becomes obvious if 
the mapping from the isotopy classes of transverse knots
in $\mathcal T$ to the isotopy classes of framed knots in $\mathcal F$ is
surjective.
However in general this mapping is not surjective and to obtain the
proof of Theorem~\ref{isomorphismtransverse} one follows the ideas of the proof of
Theorem~\ref{isomorphism}
(The famous Bennequin inequality shows that this mapping is not surjective
even for the standard contact $\R^3$.)

Thus to obtain the isomorphism between the groups $W_n^{\mathcal F}$ and
$V_n^{\mathcal T}$ it suffices to show that statement $\bf a$ of
Theorem~\ref{isomorphismtransverse} is true for the connected components 
$\mathcal T$
and $\mathcal F$ of the spaces of transverse and of framed curves. 

In~\cite{FuchsTabachnikov} D.~Fuchs and S.~Tabachnikov 
showed that statement $\bf a$ holds for 
all the connected components of the space of transverse curves 
in the case where $M$ is the standard contact $\R^3$ 
and $\mathcal A=\C$. (One can verify that the proof of this
Theorem of Fuchs and Tabachnikov goes through for $\mathcal A$ being any 
Abelian group.)
They later observed~\cite{FuchsTabachnikovprivate} that since their proof 
of this fact is mostly local, the similar
fact should be true for a big class of contact manifolds. 
However in fact the proof of their theorem is not completely local and is
based on the existence of a well-defined Bennequin invariant for a
transverse knot in $\R^3$. Unfortunately the Bennequin invariant is not 
well-defined unless the knot is zero homologous. And the generalization 
of this Theorem to the case of manifolds other than $\R^3$ meets certain
difficulties that are similar to the ones we meet when we generalize the
analogous Theorem of Fuchs and Tabachnikov for Legendrian knots,
see~\ref{reasonstoappear}. We imitate the arguments we use
in~\ref{reasonstoappear} and obtain that 
statement {\bf a} of Theorem~\ref{isomorphismtransverse}  is true for
any contact $3$-manifold with a parallelized contact structure.
(Observe that in the case of transverse knots, on the contrary to the case of
Legendrian knots, no extra conditions on the contact manifold 
appear to be needed for the statement {\bf a} to be true.)
Thus we get the following Theorem.

\begin{thm}\label{isomorphismobtainedtransverse}
Let $(M, C)$ be a contact manifold with a
parallelized contact structure, and let $\mathcal T$ be a connected component of
the space of transverse curves in $M$. Let $\mathcal F$ be the
connected component of the space of framed curves that contains $\mathcal T$.
Let $V_n^{\mathcal T}$ (resp. $W_n^{\mathcal F}$) 
be the group of $\mathcal A$-valued order $\leq n$ invariants of transverse
(resp. framed) knots from $\mathcal T$ (resp. from $\mathcal F$). 
Then the groups $V_n^{\mathcal T}$ and $W_n^{\mathcal F}$ are canonically
isomorphic.
\end{thm}

This generalizes the Theorem of D.~Fuchs and
S.~Tabachnikov~\cite{FuchsTabachnikov} saying that for the standard contact
$\R^3$
and for $\mathcal A=\C$ the quotient groups
$V_n^{\mathcal T}/V_{n-1}^{\mathcal T}$ and $W_n^{\mathcal
F}/W_{n-1}^{\mathcal F}$ are canonically isomorphic.
 
The proof of this
Theorem of Fuchs and Tabachnikov was based on the fact that for the
$\C$-valued Vassiliev invariants of framed knots in
$\R^3$ there exists the universal Vassiliev invariant
constructed by T.~Q.~T.~Le and J.~Murakami~\cite{LeMurakami}.
(For unframed knots in $\R^3$ the construction of the universal Vassiliev
invariant is the classical result of M.~Kontsevich~\cite{Kontsevich}, and
the invariant itself is the well-known Kontsevich integral.)
The existence of the universal Vassiliev invariant is currently known only
for a
very limited collection of $3$-manifolds, and only for $\mathcal A$ being 
$\C$, $\R$ or $\Q$.
(Andersen, Mattes, Reshetikhin~\cite{AMR} proved its existence in the case
where
$\mathcal A=\C$ and $M$ is the total space of an $\R^1$-bundle over a
compact oriented surface $F$ with $\p F\neq \emptyset$.)

Thus the approach used in~\cite{FuchsTabachnikov} to show the isomorphism of
the quotient groups is not applicable for almost all contact $3$-manifolds
and Abelian groups $\mathcal A$, and
Theorem~\ref{isomorphismobtainedtransverse} appears to be a strong
generalization of the result of Fuchs and Tabachnikov.

\begin{rem}
Let $(M,C)$ be a contact manifold with a parallelized contact structure, and
let $\mathcal F$ be a connected component of the space of framed curves in
$M$. Theorem~\ref{isomorphismobtainedtransverse} 
implies that for any $n\in\N$ the group of $\mathcal A$-valued order 
$\leq n$ invariants of transverse knots from a connected component 
of the space of transverse curves contained in $\mathcal F$ does not depend on the
choice of a parallelized contact structure. Hence this group can not be
used to distinguish parallelized contact structures on $M$.
\end{rem}

\section{Examples of Legendrian knots that are distinguishable by finite
order invariants.}\label{nonisomorphicgroups}
In this section we construct a big class of examples when Vassiliev
invariants distinguish Legendrian knots that realize isotopic framed knots
and are homotopic as Legendrian curves. 
Theorem~\ref{isomorphism} says that in these examples the groups of
Vassiliev invariants of Legendrian and of framed knots are not canonically 
isomorphic,
and we obtain the first known examples when these groups are not canonically
isomorphic.

Theorem of R.~Lutz~\cite{Lutz} says that for an arbitrary
orientable $3$-manifold $M$ every homotopy class of distributions of
$2$-planes tangent to $M$ contains a contact structure.
(The Theorem of Ya. Eliashberg~\cite{Eliashberg} says even more that every homotopy 
class of the distributions of $2$-planes tangent to $M$ contains a positive
overtwisted contact structure.)

However in our constructions we will use only the Euler classes of 
contact bundles. For this reason we start with the following Proposition.

\begin{prop}\label{existcontact} 
Let $M$ be an oriented $3$-manifold and let $e$ be an element of
$H^2(M,\Z)$. Then $e\in H^2(M,\Z)$ can be realized 
as the Euler class of a cooriented 
contact structure on $M$ if
and only if $e=2\alpha$, for some $\alpha \in H^2(M,\Z)$.
\end{prop}

For the Proof of Proposition~\ref{existcontact} see
Subsection~\ref{Proofexistcontact}.

\subsection{Examples of nonisotopic Legendrian knots in $S^1\times S^2$ 
that can be distinguished by Vassiliev invariants.}

Let $C$ be a cooriented contact structure on $M=S^1\times S^2$ such that
the Euler class of the contact bundle is nonzero. (The Euler class takes
values in $\Z=H^2(S^1\times S^2)$, and Proposition~\ref{existcontact} says
that for any even $i\in\Z$ there exists a cooriented contact
structure on $S^1\times S^2$ with the Euler class $i$.)

\begin{figure}[htbp]
 \begin{center}
  \epsfxsize 8cm
  \hepsffile{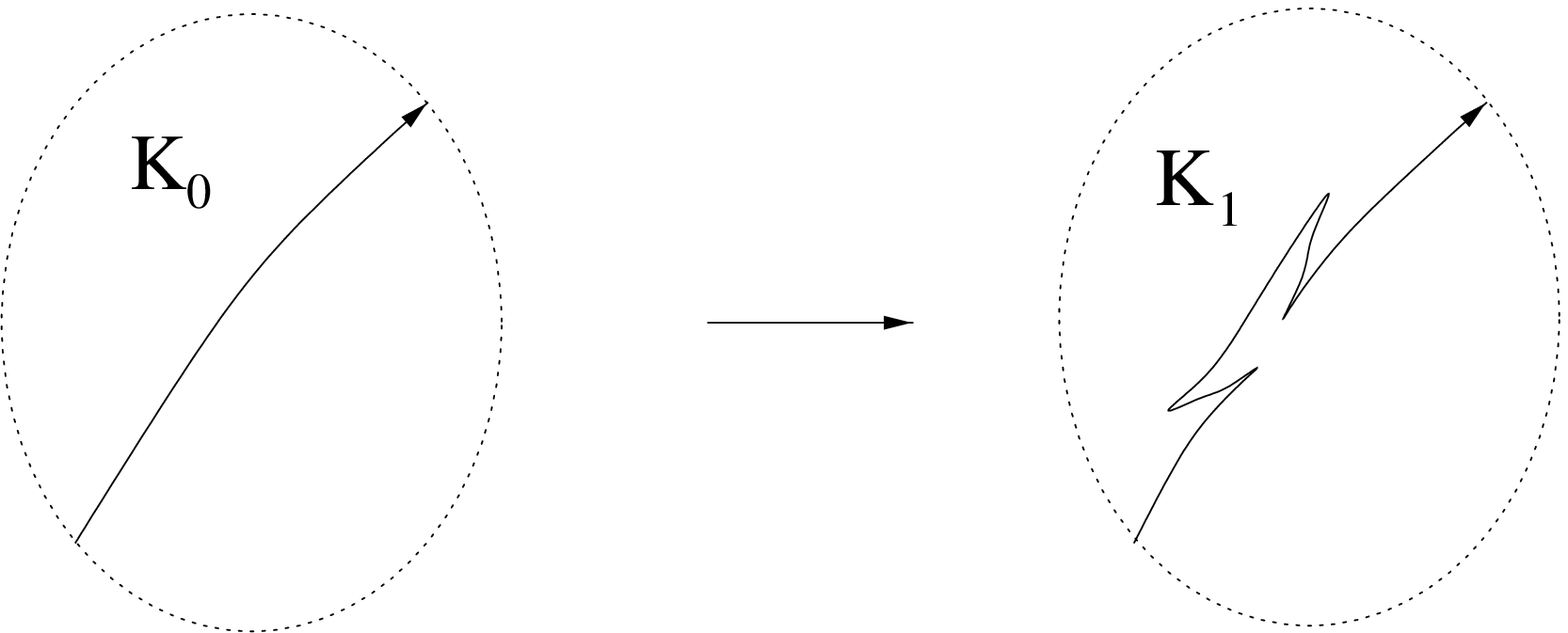}
 \end{center}
\caption{}\label{change.fig}
\end{figure}

Let $K$ be a knot in $S^1\times S^2$ that crosses exactly once one of the
spheres $t\times S^2$. The Theorem of Chow~\cite{Chow} and
Rashevskii~\cite{Rashevskii} says
that there exists a Legendrian knot $K_0$ that is $C^0$-small 
isotopic to $K$ as an unframed knot. Let $K_1$ be the Legendrian knot that
is the same as $K_0$ everywhere except of a small piece located in a chart
contactomorphic to the standard contact $\R^3$ where it is changed as it is
shown in Figure~\ref{change.fig} (see~\ref{description}).

\begin{thm}\label{example1} 
\begin{description}
\item[a] Legendrian knots $K_0$ and $K_1$ belong to the same component of
the space of Legendrian curves and realize isotopic framed knots.
\item[b] There exists a $\Z$-valued order one invariant $I$ of Legendrian
knots, such that $I(K_0)\neq I(K_1)$.
\end{description}
\end{thm}

For the Proof of Theorem~\ref{example1} see Subsection~\ref{Proofexample1}.

\begin{rem}

Let $K_i$, $i\in\N$, be the knot that is the same as $K_0$ 
everywhere except of a small piece located in a chart contactomorphic 
to the standard contact $\R^3$ where it is changed in the way described 
by the addition of $i$ zigzags shown in Figure~\ref{change.fig}. 
The Proof of Theorem~\ref{example1} implies that all $K_i$'s are homotopic
as Legendrian curves and realize isotopic framed knots, but for all $i_1\neq i_2$ 
Legendrian knots $K_{i_1}$ and $K_{i_2}$ are not Legendrian isotopic. 
The order one invariant of Legendrian knots 
$I$ constructed in the Proof of Theorem~\ref{example1}
has the property that 
$I(K_{i_1})=I(K_{i_2})+(i_2-i_1)$. Hence this $I$ distinguishes all the 
$K_i$'s.
\end{rem}

\begin{emf}\label{overtwisted}{\bf Examples of nonisotopic Legendrian knots
with overtwisted complements that realize isotopic framed knots and are
homotopic as Legendrian immersions.}
Let $\Delta$ be an embedded into $M$ disk centered at a point $p\in M$. The
Theorem of Eliashberg~\cite{Eliashberg} says that every homotopy class of
distributions of $2$-planes tangent to $M$ 
contains an overtwisted contact structure
that has $\Delta$ as the standard overtwisted disk. In the example of
Theorem~\ref{example1} we can start with an overtwisted contact structure that
has $\Delta$ as an overtwisted disk and with an unframed knot $K$ that is far away
from $\Delta$. Then since both $K_0$ and $K_1$ were constructed using a
$C^0$-small approximation of $K$, we can assume that they are
also far away from $\Delta$. And we have constructed examples of nonisotopic Legendrian
knots with overtwisted complements that realize isotopic framed knots and are
homotopic as Legendrian immersions. Previously such examples were unknown
and the Theorem of Ya. Eliashberg and M.~Fraser~\cite{EF} says that such
examples are impossible if the ambient manifold is $S^3$. 
\end{emf}

\subsection{Examples of nonisotopic Legendrian knots in the total spaces of
$S^1$-bundles over nonorientable surfaces that can be distinguished by
Vassiliev invariants.} 
\begin{emf}\label{OtherExamples}
{\em Below we describe another big family of examples where finite
order invariants distinguish Legendrian knots that realize isotopic framed
knots and are homotopic as Legendrian immersions.\/} 

Let $F$ be a nonorientable surface that can be decomposed as a connected sum of
the Klein bottle $K$ and a surface $F'\neq S^2$.
Let $M$ be an orientable manifold that admits a structure of a locally
trivial $S^1$-fibration $p:M\rightarrow F$. (For example one can take $M$ to be 
the spherical tangent bundle $STF$ of $F$.)

\begin{figure}[htbp]
 \begin{center}
  \epsfxsize 6cm
  \hepsffile{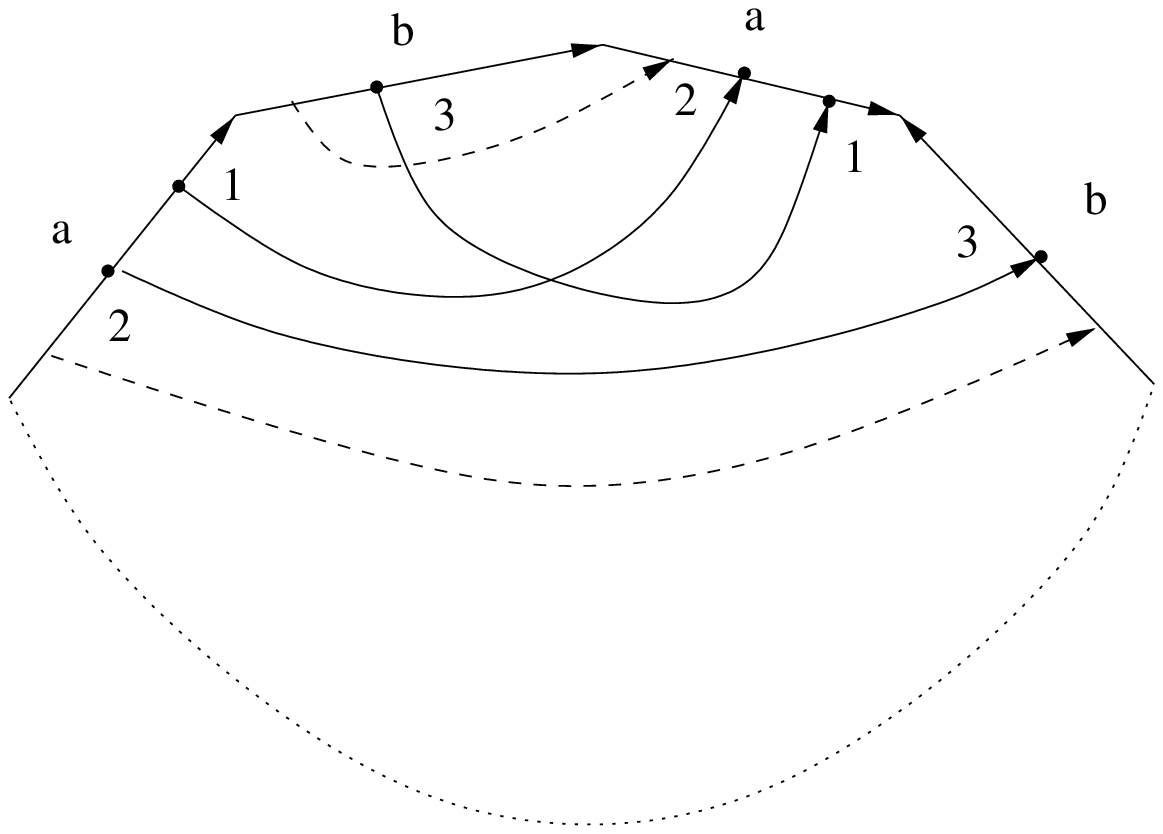}
 \end{center}
\caption{}\label{example2.fig}
\end{figure}

Consider an $S^1$-fibration $\xi:N\rightarrow S^1$ induced from $p$ by 
the mapping $S^1\rightarrow F$ that corresponds to the solid loop in
Figure~\ref{example2.fig}.
(In this Figure the enumeration of the end points
of the arcs indicates which pairs
of points should be identified to obtain the loop.)
Since the solid loop is an orientation preserving loop in $F$, we get that 
$N=T^2$ (torus). Put $\mu:N=T^2\rightarrow M$ to be the natural mapping
of the total space of the induced fibration $\xi:N\rightarrow S^1$ into the 
total space of $p:M\rightarrow F$.

A homology class in $H_1(M, \Z)$ projecting
to the dashed loop in Figure~\ref{example2.fig} has intersection $1$ 
with the class $[\mu(T^2)]\in H_2(M, \Z)$ 
realized by $\mu(T^2)$. Thus there exists $\alpha\in H^2(M,
\Z)$ such that $\alpha([\mu(T^2)])=1$.
Proposition~\ref{existcontact} says that 
for every $r\in\Z$ the class $2r\alpha$ is realizable as the Euler class
of a cooriented contact structure on $M$. Thus for every $r\in\Z$ there
exists a cooriented contact structure on $M$ such that the value of the
Euler class of the contact bundle on $[\mu(T^2)]$ is equal to $2r$.

Let $C$ be a cooriented contact structure on $M$ such that the Euler class 
$e\in H^2(M, \Z)$ of
the contact bundle satisfies $e([\mu(T^2)])=2r$, for some
nonzero $r\in\Z$.

Let $K$ be an arbitrary Legendrian knot such that its projection to 
$F$ (considered as a loop) is free homotopic to the solid loop in Figure~\ref{example2.fig}.
Let $K_1, K_2$ be Legendrian knots that are the same as $K$ everywhere except 
of a chart (contactomorphic to the standard contact $\R^3$) where $K_1$ and
$K_2$ are different from $K$ as it is described in
Figure~\ref{example4.fig}, see~\ref{description}. 
(The number of cusps in Figure~\ref{example4.fig}
is $e([\mu(T^2)])=2r\neq 0$.)

\begin{figure}[htbp]
 \begin{center}
  \epsfxsize 12cm
  \hepsffile{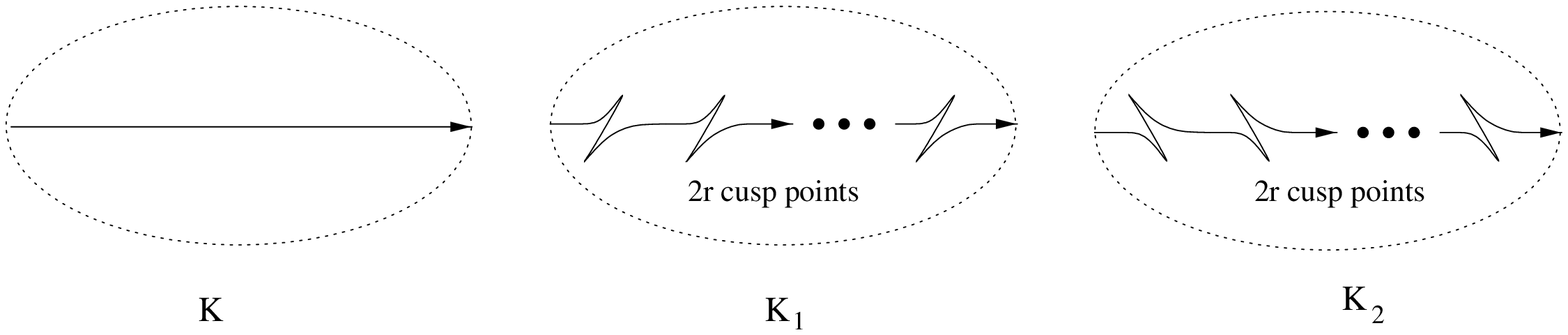}
 \end{center}
\caption{}\label{example4.fig}
\end{figure}

\end{emf}

\begin{thm}\label{example2} 
The knots $K_1$ and $K_2$ described above belong to the same
component $\mathcal L$ of the space of Legendrian curves and realize isotopic 
framed knots. There exists a $\Z$-valued order one invariant $I$ of 
Legendrian knots from $\mathcal L$ such that $I(K_1)\neq I(K_2)$.
\end{thm}

For the Proof of Theorem~\ref{example2} see Subsection~\ref{Proofexample2}.

\begin{rem} 
Similarly to~\ref{overtwisted} one verifies that the contact structure and
the knots $K_1$ and $K_2$ in the statement of Theorem~\ref{example2} can be
chosen so that the restrictions of the contact structure to the complements 
of $K_1$ and of $K_2$ are overtwisted.

Using the ideas of the Proof of Theorem~\ref{example2} one can 
construct many other examples of Legendrian knots that can be
distinguished by Vassiliev invariants of Legendrian knots
even though they realize isotopic 
framed knots and are homotopic as Legendrian immersions.
For example as 
a solid loop in Figure~\ref{example2.fig} 
we could take any loop $\beta$ such that the number of double points that
separate $\beta$ into two orientation reversing loops is odd, and the value
of the Euler class of the contact bundle on $[\mu(T^2)]\in H_2(M,\Z)$
is nonzero.
\end{rem}

\section{Proofs}
\subsection{Useful facts, Lemmas, and some technical definitions}
\begin{prop}\label{commute}
Let $p:X\rightarrow Y$ be a locally trivial $S^1$-fibration of an oriented
manifold $X$ over a (not necessarily orientable) manifold $Y$. Let
$f\in\pi_1(X)$ be the class of an oriented 
$S^1$-fiber of $p$, and let $\alpha$ be an
element of $\pi_1(X)$. Then:
\begin{description}
\item[a] $\alpha f=f\alpha\in\pi_1(X)$, provided that $p(\alpha)$ is an
orientation preserving loop in $Y$.
\item[b] $\alpha f=f^{-1}\alpha\in\pi_1(X)$, provided that $p(\alpha)$ is an
orientation reversing loop in $Y$.
\end{description}
\end{prop}

\begin{emf}{\em Proof of Proposition~\ref{commute}.\/}
If we move an oriented fiber along the loop $\alpha\in X$, then in the end it
comes to itself either with the same or with the opposite orientation. 
It is easy to see that it comes to itself with the opposite orientation if
and only if $p(\alpha)$ is an orientation reversing loop in $Y$.\qed
\end{emf}

\begin{prop}\label{Preissman}
Let $F\neq S^2, T^2\text{ (torus), } \R P^2, K\text{ (Klein bottle)}$
be a surface (not necessarily compact or orientable),
and let $G$ be a nontrivial commutative subgroup of $\pi_1(F)$.
Then $G$ is infinite cyclic.
\end{prop}

\begin{emf}{\em Proof of Proposition~\ref{Preissman}.\/}
It is well-known that any closed $F$, other than $S^2, T^2, \R P^2,K$
admits a hyperbolic metric.
(It is induced from the universal covering of $F$ by the hyperbolic plane.)
The Theorem of A.~Preissman (see~\cite{Docarmo} pp. 258-265)
says that if $M$ is a closed Riemannian manifold of negative sectional 
curvature,
then any nontrivial Abelian subgroup $G<\pi_1(M)$ is isomorphic to $\Z$.
Thus if $F\neq S^2, T^2, \R P^2, K$ is closed, then any nontrivial
commutative $G<\pi_1(F)$ is infinite cyclic.
If $F$ is not closed, then the statement of the Proposition is also
true because in this case $F$ is homotopy equivalent to a bouquet of
circles.
\qed
\end{emf}

\begin{prop}\label{toughandtechnical}
Let $F\neq S^2, \R P^2, T^2, K$ (Klein bottle) be a surface
not necessarily closed or orientable. Let $M$ be an orientable $3$-manifold, and
let $p:M\rightarrow F$ be a locally trivial $S^1$-fibration. Let $f\in\pi_1(M)$
be the class of an oriented $S^1$-fiber of $p$, and let
$\alpha\in\pi_1(M)$ be an element with $p_*(\alpha)\neq 1\in\pi_1(F)$.
Let $\beta$ be an element of the centralizer $Z(\alpha)<\pi_1(M)$ of $\alpha$.
Then
there
exist $i,j\in\Z$ and nonzero $n\in\Z$ such that $\beta^n=\alpha^if^j$.
\end{prop}

\begin{emf}\label{prooftoughandtechnical}
{\em Proof of Proposition~\ref{toughandtechnical}.\/}
Since $\alpha$ and $\beta$ commute in
$\pi_1(M)$ we get that $p_*(\alpha)$ and $p_*(\beta)$ commute in $\pi_1(F)$.
Proposition~\ref{Preissman} and the fact that $p_*(\alpha)\neq 1\in\pi_1(F)$
imply that there exist $g\in\pi_1(M)$ with
$p_*(g)\neq 1 \in\pi_1(F)$, $i\in\Z$, and nonzero $n\in\Z$
such that $p_*(g)^n=p_*(\alpha)$ and $p_*(g)^i=p_*(\beta)$.

Hence (see Proposition~\ref{commute}),
$\alpha=g^n f^k \text{ and } \beta=g^i f^l, \text{ for some } k,l\in\Z$.
Using~\ref{commute} we get that
$\beta^{n}=\alpha ^{i} f^j,\text{ for some } j\in\Z$.
Since $n$ was initially chosen to be nonzero,  we get the statement of the
Proposition. \qed
\end{emf}

\begin{emf}{\em An important homomorphism.}\label{Hansen} 
Let $X$ be a manifold, let $\Omega X$ be the space of free loops
in $X$, and let $\omega\in\Omega X$ be a loop. An element $\alpha\in\pi_1
(\Omega X,\omega)$ is realizable by a mapping 
$\mu^{\alpha}:T^2=S^1\times S^1\rightarrow X$ 
with $\mu^{\alpha}\big|_{t\times S^1}=\alpha(t)$. Let 
$t(\alpha)=\mu^{\alpha}\big|_{S^1\times 1}\in\pi_1(X, \omega(1))$ 
be the element corresponding to the trace
of the point $1\in S^1$ under the homotopy of $\omega$ described by $\alpha$.
Let $t:\pi_1(\Omega X,\omega)\rightarrow \pi_1(X, \omega(1))$ be the
homomorphism that maps $\alpha\in\pi_1(\Omega X,\omega)$ to
$t(\alpha)\in\pi_1(X, \omega(1))$.

Since the $2$-cell of $T^2$ is glued to the $1$-skeleton along
the commutation relation of the meridian and of the longitude of $T^2$,   
we get that 
$t:\pi_1(\Omega X,\omega)\rightarrow \pi_1(X, \omega(1))$ is a surjective
homomorphism of $\pi_1(\Omega X,\omega)$ onto the centralizer $Z(\omega)$
of $\omega\in\pi_1(X,\omega(1))$. 

If $t(\alpha)=t(\beta)\in\pi_1(X, \omega(1))$ for  
$\alpha, \beta\in\pi_1(\Omega X, \omega)$, then the mappings $\mu^{\alpha}$
and $\mu^{\beta}$ of $T^2$ corresponding to these loops can be deformed to
be identical on the $1$-skeleton of $T^2$. Clearly the obstruction for 
$\mu^{\alpha}$ and $\mu^{\beta}$ to be homotopic as mappings of $T^2$ 
(with the mapping of the $1$-skeleton of $T^2$ fixed under homotopy) 
is an element of $\pi_2(X)$ obtained by gluing together the boundaries of
the $2$-cells of the two tori. In particular we get the Proposition of 
V.~L.~Hansen~\cite{Hansen} saying that $t:\pi_1(\Omega X,\omega)\rightarrow Z(\omega)
<\pi_1(X,\omega(1))$ is an isomorphism, provided that $\pi_2(X)=0$.
\end{emf}

\begin{emf}\label{h-principleforcurves}
{\em $h$-principle for curves in $M$.}
For a $3$-dimensional manifold $M$ we put $STM$ to be the manifold obtained by 
the fiberwise spherization of the tangent bundle of $M$, and we put
$\pr':STM\rightarrow M$ to be the corresponding locally trivial $S^2$-fibration.
The $h$-principle (that can be found in~\cite{Gromov})
says that the space of curves in $M$ is weak homotopy
equivalent to $\Omega STM$ (the space of free loops in $STM$). The weak
homotopy equivalence is given by mapping a curve $K$ to a loop $\vec
K\in\Omega STM$ that sends a point $t\in S^1$ to the point of $STM$
corresponding to the direction of the velocity vector of $K$ at $K(t)$.
\end{emf}

\begin{defin}[of $m(K_1, K_2)$ and of $K^0, K^{\pm 1}, K^{\pm
2}\dots$]\label{obstruction} 
Let $K_1$ and $K_2$ be two framed knots that 
coincide pointwise as embeddings of $S^1$. Then there is an integer obstruction 
$m(K_1, K_2)\in\Z$ for them to be isotopic as framed knots with the embeddings 
of $S^1$ fixed under the isotopy. This obstruction is calculated as follows.
Let $K'_1$ be the knot obtained by shifting $K_1$ along the framing and
reversing the orientation on the shifted copy. Together $K_1$ and $K'_1$
bound a thin strip. We put $m(K_1, K_2)$ to be the intersection number of 
the strip with a very small shift of $K_2$ along its framing. 

For a framed knot $K^0$ we denote by $K^i$, $i\in \Z$, the isotopy class of a
framed knot that coincides with $K^0$ as an embedding of $S^1$ and has
$m(K^0, K^i)=i$.

For two singular framed knots $K_{1s}$ and $K_{2s}$ with $n$ transverse 
double points
that coincide pointwise as immersions of $S^1$, we put $m(K_{1s},
K_{2s})\in\Z$ to be the value of $m$ on the nonsingular framed knots $K_1$ and $K_2$
that coincide pointwise as embeddings of $S^1$ and are obtained from
$K_{1s}$ and $K_{2s}$ by resolving each pair of the corresponding double points
of $K_{1s}$ and of $K_{2s}$ in the same way. (The value $m(K_{1s}, K_{2s})$
does not depend on the resolution as soon as the corresponding double points
of the two knots are resolved in exactly the same way.)
As before $m(K_{1s}, K_{2s})$ is the integer valued obstruction for $K_{1s}$
and $K_{2s}$ to be isotopic as singular framed knots 
with the immersion of $S^1$ corresponding to the two knots fixed under isotopy.

For a  singular framed knot $K^0_s$ with $n$ transverse 
double  points we denote by $K^i_s$, $i\in\Z$,
the isotopy class of a singular framed knot with $n$ transverse double points that
coincides with $K^0_s$ as an immersion of $S^1$ and has $m(K_s^0, K_s^i)=i$.
\end{defin}

\begin{prop}\label{homotopy} Let $K_1$ and $K_2$ be framed knots 
(resp. singular framed knots with $n$ transverse double points) that coincide
pointwise as embeddings (resp. immersions) of $S^1$. 
Then $K_1$ and $K_2$ are homotopic as
framed knots (resp. singular framed knots with $n$ transverse double points)
if and only if $m(K_1, K_2)$ is even.
\end{prop}

\begin{emf}\label{proofhomotopy}{\em Proof of Proposition~\ref{homotopy}.\/}
Clearly if $m(K_1, K_2)$ is even, then $K_1$ and $K_2$ are framed
homotopic. (We can change the obstruction by two by creating a small kink
and passing through a double point at its vertex.)

Every oriented $3$-dimensional manifold $M$ is parallelizable, and hence it 
admits a $\spin$-structure. A framed
curve $K$ in $M$ represents a loop in the principal $SO(3)$-bundle of $TM$. 
(The $3$-frame corresponding to a point of $K$ is the velocity
vector, the framing vector, and the unique third vector of unit length such
that the $3$-frame defines the positive orientation of $M$.)
One
observes that the values of the $\spin$-structure on the loops in the principal
$SO(3)$-bundle of $TM$ realized by $K_1$ and $K_2$ are different provided that 
$m(K_1, K_2)$ is odd. But these values do not change under homotopy of
framed curves. Hence if $m(K_1, K_2)$ is odd, then $K_1$ and $K_2$ are not
framed homotopic.
\qed
\end{emf} 
 
\begin{defin}[of the number of framings of a knot]\label{strange} 
Using the self-linking invariant of framed knots one can easily
show that if $K_1$ and $K_2$ in~\ref{obstruction} are pointwise coinciding 
zero-homologous framed knots and  $m(K_1, K_2)\neq 0$,
then $K_1$ is not isotopic to $K_2$ in the category of framed knots. However
for knots that are not zero-homologous 
this is not generally true, see~\ref{Proofaexample1}.
For this reason we introduce the following definitions.

If for an unframed knot $K$ there exist isotopic framed knots $K_1$ and $K_2$ 
that coincide with $K$ pointwise and have $m(K_1, K_2)\neq 0$, then we say that 
$K$ {\em admits finitely many framings\/}. 
For $K$ that admits finitely many
framings we put {\em the number of framings $m_K$ of\/} 
$K$ to be the minimal positive
integer $l$ such that there exist isotopic framed knots $K_1$ and $K_2$ that 
coincide with $K$ pointwise and have $m(K_1, K_2)=l$.
One can easily show that if $K$ admits finitely many framings, then there are
exactly $m_K$ isotopy classes of framed knots realizing the isotopy class of
the unframed knot $K$. Proposition~\ref{homotopy} implies that $m_K$ is even.

In a similar way we introduce the notion of the number of framings 
for unframed singular knots with $n$ double points.
\end{defin}

\begin{prop}\label{decrease} 
Let $(M,C)$ be a contact $3$-manifold with a cooriented contact structure,
let $\mathcal F$ be a connected component of the space of framed curves, and  
let $\mathcal L\subset \mathcal F$ be a connected component of the space of Legendrian curves
in $(M, C)$.
\begin{description}
\item[a] Let $K$ be an unframed knot obtained by forgetting the framing on a
knot from $\mathcal F$. Then there exists a Legendrian knot from $\mathcal
L$ realizing the isotopy class of $K$.
\item[b] If $K^0$ is an isotopy class of framed knots in $\mathcal F$ that 
is realizable by a Legendrian knot from $\mathcal L$, then the isotopy 
class of $K^{-2}$ (see~\ref{obstruction}) is also realizable by a Legendrian knot from $\mathcal L$.
\item[c] Let $K_s$ be an unframed singular knot with $n$ double  points 
obtained by forgetting the framing on a singular knot  
from $\mathcal F$. Then there exists a singular Legendrian knot from $\mathcal  
L$ realizing the isotopy class of $K_s$.
\item[d] If $K_s^0$ is an isotopy class of singular 
framed knots in $\mathcal F$ that  is realizable by a singular 
Legendrian knot from $\mathcal L$, then the isotopy class of $K_s^{-2}$ is
also realizable by a singular Legendrian knot from $\mathcal L$.
\end{description}
\end{prop}

\begin{emf}\label{proofdecrease}
{\em Proof of statement {\bf a} of Proposition~\ref{decrease}.\/}
Let $CM$ be the fiberwise spherization of the $2$-dimensional
contact vector bundle, and let $\pr:CM\rightarrow M$ be the corresponding 
locally trivial $S^1$-fibration. We denote by $f\in \pi_1(CM)$ the class of
an oriented $S^1$-fiber of $\pr$. 
For a Legendrian curve
$K_l:S^1\rightarrow M$ denote by $\vec K_l$ the loop in $CM$ obtained by
mapping a point $t\in S^1$ to the point of $CM$ corresponding to the
direction of the velocity vector of $K_l$ at $K_l(t)$. 

The $h$-principle~\ref{h-principleLegendrian} 
says that Legendrian
curves  $K_1$ and $K_2$ in $M$ belong to the same component of the space of Legendrian
curves in $M$ if and only if $\vec K_1$ and $\vec K_2$ are
free homotopic loops in $CM$. 

W.~L. Chow~\cite{Chow} and P.~K.~Rashevskii~\cite{Rashevskii} showed that
every unframed knot $K$ is isotopic to a Legendrian knot $K_l$
(and this isotopy can be made $C^0$-small). 
Deforming $K$ we can assume (see~\ref{h-principleforcurves}) that: {\bf 1:} $K$ and $K_l$ coincide in
the neighborhood of $1\in S^1\subset \C$, {\bf 2:} $K$ and $K_l$ realize the
same element $[K]\in \pi_1(M, K_l(1))$, and {\bf 3:} that liftings to $CM$ of 
Legendrian curves from $\mathcal L$ are free homotopic to a loop $\alpha$ in
$CM$ such that $\alpha(1)=\vec
K_l(1)$ and $\pr (\alpha)=[K] \in \pi_1(M, K_l(1))$. 

Proposition~\ref{commute} says that $f$ is in the
center of $\pi_1(CM, \vec K_l(1))$, since the contact structure is
cooriented and hence oriented. Then $\vec K_l=\alpha f^i\in \pi_1(CM, \vec K_l
(1))$, for some $i\in\Z$. 

Take a chart of $M$ (that is contactomorphic to the standard contact $\R^3$) 
containing a piece of the Legendrian knot.
From the formula for the Maslov number deduced 
in~\cite{FuchsTabachnikov} it is easy to see that the modifications of the
Legendrian knot corresponding to the insertions of 
two cusps shown in Figure~\ref{twocusp.fig} (see~\ref{description})
induce multiplication by $f^{\pm 1}$ of the lifting of $K_l$ to an
element of $\pi_1(CM,\vec K_l(1))$. (Here the sign depends on
the choice of an orientation of the fiber used to define $f$.) Performing
this operation sufficiently many times we obtain the Legendrian knot from
$\mathcal L$ realizing the isotopy class of the unframed knot $K$.

\begin{figure}[htbp]
 \begin{center}
  \epsfxsize 10cm
  \hepsffile{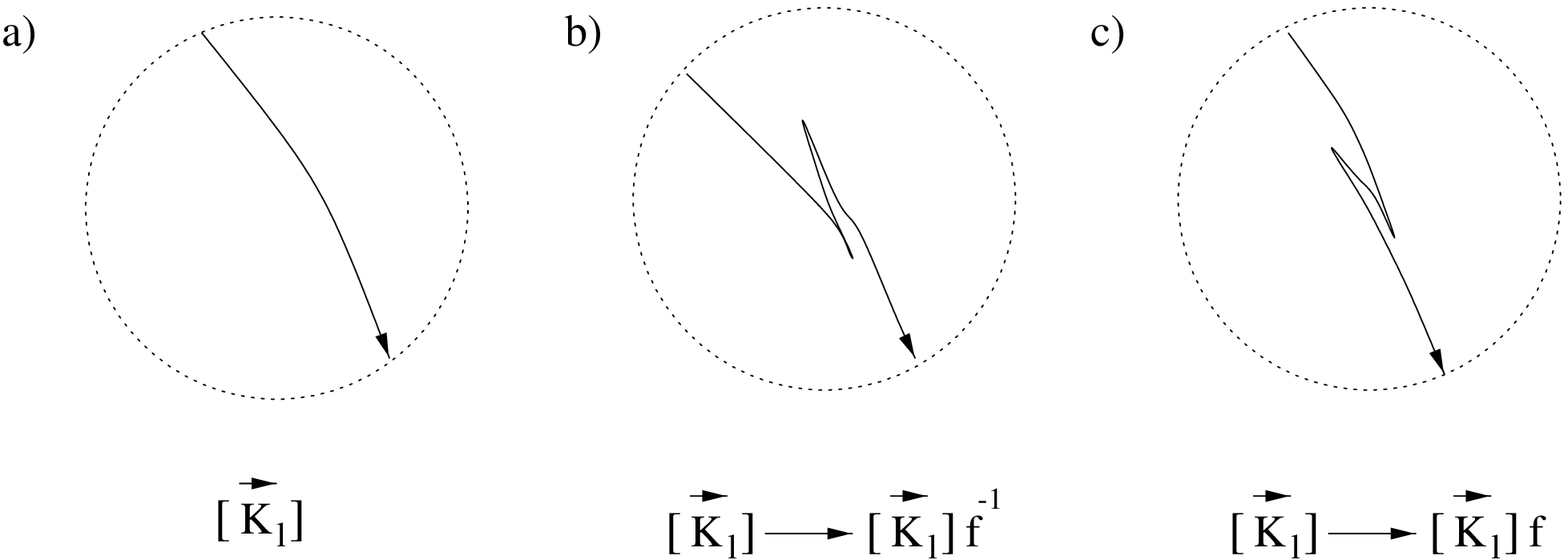}
 \end{center}
\caption{}\label{twocusp.fig}
\end{figure}

One easily modifies the arguments above to obtain 
the proof of statement {\bf c} of Proposition~\ref{decrease}.

{\em Proof of statement {\bf b} of Proposition~\ref{decrease}.\/}
Take a chart of $M$ (that is contactomorphic to the standard contact
$\R^3$) containing a piece of the knot $K^0$ and perform the
homotopy in $\mathcal L$ shown in Figure~\ref{kink.fig},
see~\ref{description}. (Observe that a
self-tangency point of the projection of a Legendrian curve in $\R^3$ to the
$(x,z)$-plane corresponds to a double point of the Legendrian curve.)
Straightforward verification (cf. the formula for the Bennequin invariant
deduced in~\cite{FuchsTabachnikov}) shows that the Legendrian knot we 
obtain in the end of the homotopy realizes $K^{-2}$. 

\begin{figure}[htbp]
 \begin{center}
  \epsfxsize 12cm
  \hepsffile{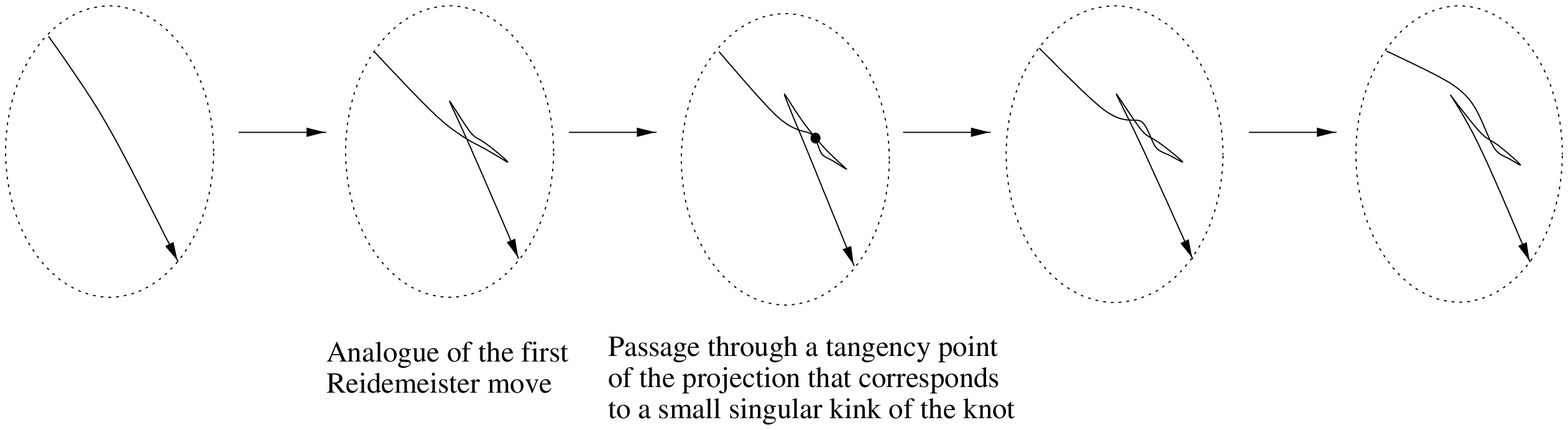}
 \end{center}
\caption{}\label{kink.fig}
\end{figure}

One easily modifies these arguments to obtain the proof of statement {\bf d} of
Proposition~\ref{decrease}.
This  finishes the proof of Proposition~\ref{decrease}.\qed 
\end{emf}

\subsection{Proof of Theorem~\ref{isomorphism}}\label{proofisomorphism}
The fact that statement {\bf b} of Theorem~\ref{isomorphism} implies 
statement {\bf a} is clear.
Thus we have to  show that statement {\bf a} implies statement
\textrm{\bf b}. This is done by showing that there exists a homomorphism 
$\psi:V_n^{\mathcal L}\rightarrow W_n^{\mathcal F}$ such that 
$\phi\circ\psi=\id_{V_n^{\mathcal L}}$ and $\psi\circ\phi=\id_{W_n^{\mathcal
F}}$.

Let $x\in V_n^{\mathcal L}$ be an invariant.
In order to construct $\psi(x)\in W_n^{\mathcal F}$ we have to
specify the value of $\psi (x)$ on every framed knot $K\in\mathcal F$.

\begin{emf}\label{definitionofpsi}{\em Definition of $\psi(x)$.\/} If the 
isotopy class of the knot $K\in \mathcal
F$ is realizable by a Legendrian knot $K_l\in\mathcal L$, then put
$\psi(x)(K)=x(K_l)$. The value $\psi(x)(K)$ is well-defined because if
$K'_l\in\mathcal L$ is another knot realizing $K$, then
$x(K_l)=x(K'_l)$ by statement {\bf a} of Theorem~\ref{isomorphism}.

Let $\mathcal C$ be the component of the space of unframed curves
that corresponds to forgetting framings on the curves from $\mathcal F$.
Propositions~\ref{decrease} and~\ref{homotopy} imply that 
if an unframed knot $K_u\in \mathcal C$
admits finitely many framings (see~\ref{strange}), 
then all the isotopy classes of framed knots
from $\mathcal F$ realizing the isotopy class of the unframed knot $K_u$ are realizable by
Legendrian knots from $\mathcal L$. Thus we have defined the value of
$\psi(x)$ on all the framed knots from $\mathcal F$ that realize unframed 
knots admitting finitely many framings. 

If $K_u\in \mathcal C$ admits infinitely many framings, then either {\bf 1)} 
all the isotopy classes of 
framed knots from $\mathcal F$ realizing the isotopy class of $K_u$ are
realizable by Legendrian knots from $\mathcal L$ or {\bf 2)} there exists a
knot $K^0\in \mathcal F$  realizing the isotopy class of $K_u$ 
such that $K^0$ is realizable by a Legendrian knot
from $\mathcal L$ and $K^{+2}$ (see~\ref{obstruction})
is not realizable by a Legendrian knot from 
$\mathcal L$. (In this case $K^{+4}, K^{+6}$ etc. also are not realizable by
Legendrian knots from $\mathcal L$, see~\ref{decrease}.) In the case {\bf 1)} 
the value of $\psi(x)$ is already defined on 
all the framed knots from $\mathcal F$ realizing $K_u$. In the case {\bf 2)} put 
\begin{equation}\label{eqextension}
\psi(x)(K^{+2})=\sum_{i=1}^{n+1}\bigl((-1)^{i+1}\frac{(n+1)!}{i!(n+1-i)!}
\psi(x)(K^{+2-2i})\bigr).
\end{equation}
(Proposition~\ref{decrease} implies that the sum on the right hand side 
is well-defined.)
Similarly put 
\[
\begin{array}{l}

\psi(x)(K^{+4})=\sum_{i=1}^{n+1}\bigl((-1)^{i+1}\frac{(n+1)!}{i!(n+1-i)!}
\psi(x)(K^{+4-2i})\bigr), \\
\psi(x)(K^{+6})=\sum_{i=1}^{n+1}\bigl((-1)^{i+1}\frac{(n+1)!}{i!(n+1-i)!}
\psi(x)(K^{+6-2i})\bigr)\text{ etc.}\\
\end{array}
\]

Now we have defined $\psi(x)$ on all the framed knots (from $\mathcal F$) realizing $K_u$.
Doing this for all $K_u$ for which case {\bf 2)} holds we define the value of $\psi(x)$
on all the knots from $\mathcal F$.

\end{emf}

{\bf Below we show that $\psi(x)$ is an order $\leq n$ invariant of framed knots from
$\mathcal F$.}
We start by proving the following Proposition.

\begin{prop}\label{mainidentity}
Let $K^0$ be a framed knot from $\mathcal F$, then $\psi(x)$ defined as
above satisfies identity~\eqref{eqextension}.
\end{prop}

\begin{emf}{\em Proof of Proposition~\ref{mainidentity}.\/}
If $K^{+2}$ is not realizable by a Legendrian knot from $\mathcal L$, then
the statement of the proposition follows from the formula we used to
define $\psi(x)(K^{+2})$. 

If $K^{+2}$ is realizable by a Legendrian knot $K_l^0$, then consider a
singular Legendrian knot $K_{ls}$ with $(n+1)$ double points that are
vertices of $(n+1)$ small kinks such that we get $K_l$
if we resolve all the double points positively staying in the class of the Legendrian
knots. (To create $K_{ls}$ we perform the first half of the
homotopy shown in Figure~\ref{kink.fig} in $n+1$ places on $K_l^0$.)

Let $\Sigma$ be the set of the $2^{n+1}$ possible resolutions of the double
points of $K_{ls}$. 
For $\sigma\in \Sigma$ put $\sign(\sigma)$ to be the sign of the resolution, and
put $K_{ls}^{\sigma}$ to be the nonsingular Legendrian knot obtained via the
resolution $\sigma$. Since $x$ is an order $\leq n$ invariant of Legendrian 
knots we get that
\begin{equation}
0=\sum_{\sigma\in \Sigma}\bigl(\sign(\sigma)x(K_{ls}^{\sigma})\bigr)=
\psi(x)(K_l^0)+\sum_{i=1}^{n+1}(-1)^i\frac{(n+1)!}{i!(n+1-i)!}\psi(x)(K_l^{-2i}).
\end{equation}
(Observe that if we resolve $i$ double points of $K_{ls}$ negatively, then
we get the isotopy class of $K_l^{-2i}$.)
This finishes the proof of the Proposition. \qed
\end{emf}

\begin{emf}
Let $K_s\in\mathcal F$ be a singular framed knot with $(n+1)$ double points. Let $\Sigma$
be the set of the $2^{n+1}$ possible resolutions of the double
points of $K_s$. For $\sigma\in
\Sigma$ put $\sign(\sigma)$ to be the sign of the resolution, and put $K^{\sigma}_s$
to be the isotopy class of the knot obtained via the resolution $\sigma$.

In order to prove that $\psi(x)$ is an order $\leq n$ invariant
of framed knots from $\mathcal F$, we have to show that 
\begin{equation}\label{identitytoshow}
0=\sum_{\sigma\in \Sigma}\bigl(\sign(\sigma)\psi(x)(K_{s}^{\sigma})\bigr),
\end{equation}
for every $K_s\in \mathcal F$.

First we observe that the fact whether identity~\eqref{identitytoshow}
holds or not depends only on the isotopy class of the singular knot $K_s$ 
with $(n+1)$ double points. 

If the isotopy class of $K_s$ is realizable by 
a singular Legendrian knot from $\mathcal L$, then
identity~\eqref{identitytoshow} holds for $K_s$, since $x$ is an order $\leq n$
invariant of Legendrian knots (and the value of $\psi(x)$ on a framed knot
$K\in\mathcal F$ realizable by a Legendrian knot $K_l\in\mathcal L$ was put to be $x(K_l)$).

Proposition~\ref{decrease} says that the isotopy class of the 
singular unframed knot 
$K_{us}$ obtained by forgetting the framing on $K_s$ is realizable by a 
singular Legendrian knot from $\mathcal L$. 

If $K_{us}$ admits finitely many framings, 
then all the isotopy classes of singular framed knots from $\mathcal F$ realizing the isotopy
class of $K_{us}$ are realizable by singular 
Legendrian knots from $\mathcal L$, and we get that
identity~\eqref{identitytoshow} holds for $K_s$. 

If $K_{us}$ admits infinitely 
many framings and all the isotopy classes of 
singular framed knots from $\mathcal F$ realizing $K_{us}$ 
are realizable by singular Legendrian knots from $\mathcal L$, 
then~\eqref{identitytoshow} automatically holds for $K_s$. 
If $K_{us}$ admits infinitely many framings 
but not all the isotopy classes of singular framed knots from $\mathcal F$ 
realizing $K_{us}$ are 
realizable by singular Legendrian knots from $\mathcal L$. Then put $K_{us}^0$ to be the framed
knot realizing $K_{us}$ that is realizable by a singular Legendrian knot
from $\mathcal L$
and such that 
$K_{us}^{2i}$, $i>0$, are not realizable by singular Legendrian knots from
$\mathcal L$. 

Proposition~\ref{decrease} says that $K_{us}^{-2i}$, $i>0$, are realizable by
singular Legendrian knots from $\mathcal L$ and hence 
identity~\eqref{identitytoshow} holds for $K_{us}^{-2i}$, $i\geq 0$. 
Using Proposition~\ref{mainidentity} and 
the fact that identity~\eqref{identitytoshow} holds for
$K_{us}^{-2i}$, $i\geq 0$, we 
show that~\eqref{identitytoshow} holds for $K_{us}^{+2}$. Namely, 
\begin{multline}
\sum_{\sigma\in
\Sigma}\sign(\sigma)\psi(x)({K_{us}^{+2}}^{\sigma})\\=
\sum_{\sigma\in
\Sigma}\Bigl(\sign(\sigma)
\sum_{i=1}^{n+1}(-1)^{i+1}\frac{(n+1)!}{i!(n+1-i)!}
\psi(x)({K_{us}^{(+2-2i)}}^{\sigma})\Bigr)\\ 
=
\sum_{i=1}^{n+1}\Bigl((-1)^{i+1}\frac{(n+1)!}{i!(n+1-i)!}\times\Bigl(\sum_{\sigma\in
\Sigma}\sign(\sigma)\psi(x)({K_{us}^{(+2-2i)}}^{\sigma})\Bigr)\Bigr)\\ 
=
\sum_{i=1}^{n+1}\Bigl( (-1)^{i+1}\frac{(n+1)!}{i!(n+1-i)!}\times\Bigl(0\Bigr)\Bigr)=0.
\end{multline}

Similarly we show that ~\eqref{identitytoshow} holds for $K_{us}^{+4},
K_{us}^{+6}, \text{ etc}\dots$. 
\end{emf}

\begin{emf} Clearly $\phi\circ\psi=\id_{V_n^{\mathcal L}}$.

Considering the values of $y\in W_n^{\mathcal F}$ on the $2^{n+1}$ possible
resolutions of a singular framed knot with $n+1$ singular small kinks we get
that $y$ should satisfy identity~\eqref{eqextension}.  
Hence $\psi\circ\phi=\id_{W_n^{\mathcal
F}}$ and this finishes the proof of Theorem~\ref{isomorphism}. \qed


\end{emf}

\subsection{The reasons for the condition~$(*)$ to
appear}\label{reasonstoappear}
The proof of the Theorem of Fuchs and Tabachnikov that says that statement
$\bf a$ of Theorem~\ref{isomorphism} is true for all the connected
components of the space of Legendrian curves when the ambient contact
manifold is the standard contact $\R^3$ is
based on the following three observations: 
\begin{description}
\item[1] There are two types of cusps arising under the projection of the
part of a Legendrian knot that is contained in a Darboux chart to the $(x,z)$-plane
(see~\ref{description}). They are formed by cusps for which the 
branch of the projection of the knot going away from the cusp is locally
located respectively above or below the tangent line at the cusp point, see
Figure~\ref{twocusp.fig}. For a Legendrian knot $K$ and $i,j\in\N$ we denote by 
$K^{i,j}$ the Legendrian knot
obtained from $K$ by the modification corresponding to an addition of 
$i$ cusp pairs of the first type and $j$ cusp pairs of the second type to
the projection of the part of $K$ located in a Darboux chart.

Let $K_1$ and $K_2$ be Legendrian knots in
the standard contact $\R^3$ that realize isotopic unframed knots. Then for
any $n_1$ and $n_2$ large enough there exist $n_3, n_4\in \N$ such that 
the Legendrian knot $K_1^{n_1,n_2}$ is Legendrian isotopic to
$K_2^{n_3,n_4}$.
\item[2] If there exists $n\in\N$ such that Legendrian knots $K_1^{n,n}$ and 
$K_2^{n,n}$ are Legendrian isotopic, then every Vassiliev invariant of Legendrian 
knots takes equal values on $K_1$ and on $K_2$.
\item[3] The number $n$ from the previous observation exists if the ambient
contact manifold is $\R^3$ and the Legendrian knots $K_1$ and $K_2$ belong to the
same component of the space of Legendrian curves and realize isotopic framed
knots. 
\end{description}

The first two observation are true for any contact $3$-manifold (since the proof
of the corresponding facts is local). But the number $n$ from the 
statement
of the third observation does not exist in general. 
In the case of the
ambient manifold being $\R^3$ Fuchs and Tabachnikov showed the existence of
such $n$ using the explicit calculation involving the Maslov classes and
Bennequin invariants of Legendrian knots. However in order for the 
Bennequin invariant to be well-defined the
knots have to be zero-homologous, and in order for the Maslov class to be
well-defined the knots have to be zero-homologous or the contact structure
has to be parallelizable.

Below we show that such $n$ exists for any
$K_1$ and $K_2$ that realize isotopic framed knots and belong to the same
component of the space of Legendrian curves, provided that the connected
component $\mathcal F$ of the space of framed curves that contains $K_1$ and $K_2$
satisfies condition~$(*)$. (We assume that the contact
structure on $M$ is cooriented.)

Let $K_1$ and $K_2$ be Legendrian knots as above, and 
let $n_1, n_2, n_3, n_4\in \N$ be such that $K_1^{n_1, n_2}$ and
$K_2^{n_3,n_4}$ are Legendrian isotopic.

We start by showing that $n_1, n_2, n_3, n_4$ can be chosen so that
$n_1+n_2=n_3+n_4$, and that
if $\mathcal F$ satisfies condition~$(*)$, then 
$n_1-n_2=n_3-n_4$.

\begin{emf}{\em Proof of the fact that $n_1, n_2, n_3, n_4$ can be chosen so
that $n_1+n_2=n_3+n_4$.} 
%
%
Let $\mu:S^1\times [0,1]\rightarrow M$ be the isotopy changing $K_1$ to $K_2$ in the
category of framed knots.
Analyzing the
proof of Fuchs and Tabachnikov one verifies that for $n_1, n_2$ large enough
the Legendrian isotopy $\bar \mu$ 
changing $K_1^{n_1, n_2}$ to $K_2^{n_3, n_4}$ can be
chosen so that for every $t\in [0,1]$ the Legendrian knot $\bar\mu
_t:S^1\times t\rightarrow M$ is contained in a thin tubular neighborhood $T_t$ of
$\mu _t:S^1\times t \rightarrow M$ and is isotopic (as an
unframed knot) to $\mu _t$ inside $T_t$. 

For two framed knots $\mu _t$ and  $\bar\mu _t$ realizing unframed knots
that are isotopic inside $T_t$ there is a well-defined $\Z$-valued
obstruction to be isotopic inside $T_t$ in the category of framed knots.
This obstruction is the difference of the self-linking numbers of the
inclusions of $\mu _t$ and $\bar \mu _t$ into $\R^3$ induced by an
identification of $T_t$ with the standard solid torus in $\R^3$. (One 
verifies that for  $\mu _t$ and $\bar \mu _t$ that are isotopic as unframed
knots inside $T_t$ this difference does not depend on the choice of the
identification of $T_t$ with the standard solid torus in $\R^3$.)

From the formula for the Bennequin invariant stated in~\cite{FuchsTabachnikov} 
one gets that the value of the obstruction for $K_1^{n_1, n_2}$ to be
isotopic as a framed knot to $K_1$ inside $T_0$ is equal to $n_1+n_2$.
Similarly the value of the obstruction for $K_2^{n_3, n_4}$ to be
isotopic as a framed knot to $K_2$ inside $T_1$ is equal to $n_3+n_4$.
Clearly the value of the obstruction for $\bar \mu_t$ to be isotopic to 
$\mu_t$ inside $T_t$ does not depend
on $t$ (for the isotopy $\mu$ changing $K_1$ to $K_2$ in the category of
framed knots), and  we get that $n_1+n_2=n_3+n_4$.  
\end{emf}

\begin{emf}{\em Proof of the fact that if $\mathcal F$ satisfies
condition~$(*)$, then $n_1-n_2=n_3-n_4$.\/}
Let $f\in\pi_1(CM)$ be the class of an oriented $S^1$-fiber of
$\pr:CM\rightarrow M$. From the $h$-principles for Legendrian and for
unframed curves (see~\ref{h-principleLegendrian} and~\ref{h-principleforcurves}) 
one obtains that every component of the space of Legendrian
curves contained in $\mathcal F$ corresponds to the conjugacy class of 
$\vec K_1 f^l\in\pi_1(CM)$, for some $l\in\Z$. (Connected components of the
space of free loops in $CM$ are naturally identified
with the conjugacy classes of the elements of $\pi_1(CM)$.) 

Using Proposition~\ref{commute}
one verifies that if $\mathcal F$ satisfies condition~$(*)$, then
for every nonzero $l\in\Z$ the elements $\vec K_1$ and $\vec K_1 f^l$ are not conjugate in
$\pi_1(CM)$.
   
From the formula for the Maslov number deduced 
in~\cite{FuchsTabachnikov} and the $h$-principle for Legendrian curves one 
gets that $K_1^{n_1,n_2}$ is contained in the component of the space of Legendrian curves
that corresponds to the conjugacy class of $\vec K_1
f^{n_1-n_2}\in\pi_1(CM)$. Using the fact that $\vec K_1$ and $\vec K_2$ are
conjugate in $\pi_1(CM)$ (since $K_1$ and $K_2$ are Legendrian homotopic)
and the fact that since the contact structure is cooriented $f$ is in the center of
$\pi_1(CM)$ (see~\ref{commute}), we get that $K_2^{n_3,n_4}$ is contained in the
component that corresponds to the conjugacy class of $\vec K_1
f^{n_3-n_4}\in\pi_1(CM)$. Since $K_1^{n_1,n_2}$ and $K_2^{n_3,n_4}$ 
are Legendrian isotopic
(and hence Legendrian homotopic) we get that $\vec K_1 f^{n_1-n_2}$ and 
$\vec K_1 f^{n_3-n_4}$ are conjugate in $\pi_1(CM)$, and
using~\ref{commute} we get that $\vec K_1$ is conjugate to $\vec
K_1f^{(n_1-n_2)-(n_3-n_4)}$. But 
since $\mathcal F$ satisfies condition~$(*)$ we have $(n_1-n_2)-(n_3-n_4)=0$, and hence $n_1-n_2=n_3-n_4$. 
\end{emf}

From the identities $n_1+n_2=n_3+n_4$ and $n_1-n_2=n_3-n_4$ one gets that 
$n_1=n_3$ and $n_2=n_4$. Assume that $n_1\geq n_2$. (The case where 
$n_2>n_1$ is treated similarly.) Put $k=n_1-n_2$.
It is easy to show that
since $K_1^{n_1, n_2}$ and $K_2^{n_3,n_4}$ are Legendrian isotopic, 
then $K_1^{n_1, n_2+k}$ and $K_2^{n_3,n_4+k}$ are also Legendrian isotopic.
(Basically one can keep the $k$ extra cusp pairs
close together on a small piece of the projection of the part of 
the knot contained in a Darboux chart during the whole isotopy
process.) But $K_1^{n_1, n_2+k}$ and $K_2^{n_3,n_4+k}$ are obtained from $K_1$ and $K_2$ 
by the
modification corresponding to the addition of $n_1=n_2+k=n_3=n_4+k$ pairs of
cusps of each of the two types, and we can take $n$ from the
observation {\bf 2} to be $n_1=n_2+k=n_3=n_4+k$. 

This shows that $K_1$ and
$K_2$ can not be distinguished by the Vassiliev invariants of Legendrian
knots provided that $\mathcal F$ satisfies condition~$(*)$, and that
$K_1$ and $K_2$ realize isotopic framed knots and are homotopic as
Legendrian immersions. Hence statement {\bf a} of Theorem~\ref{isomorphism}
is true provided that $\mathcal F$ satisfies condition~$(*)$.

\subsection{Proof of Proposition~\ref{interpretationconditionII}.\/}
\label{ProofinterpretationconditionII}
The $h$-principle for curves~\ref{h-principleforcurves} says
that the set $\mathcal C$ 
of the connected components of the space of curves in $M$ is
naturally identified with the set of the connected components of the space of
free loops in the spherical tangent bundle  $STM$ of $M$. Hence it is also
naturally identified 
with the set of conjugacy classes of the elements of $\pi_1(STM)$. (From
the long homotopy sequence of the fibration $\pr':STM\rightarrow M$ we see
that it is also naturally identified with the set of conjugacy classes of
the elements of $\pi_1(M)$.) Choose a $\spin$-structure on $M$. It is easy to
see (cf.~\ref{homotopy} and~\ref{proofhomotopy}) that the set 
$\mathcal C_{\mathcal F}$ of the connected components of the space of framed curves 
in $M$ is identified with the product $\Z_2\times \mathcal C$.
Here the $\Z_2$-factor is the value of the $\spin$-structure on the loop in
the principal $SO(3)$-bundle of $TM$ that corresponds to a framed curve
from the connected component, see~\ref{proofhomotopy}. 
(This value does not depend on the choice of the framed curve in the component.)

The $h$-principle for the Legendrian curves says that the set of the
connected components of the space of Legendrian curves is naturally
identified with the set of homotopy classes of free loops in $CM$ (the
spherical contact bundle of $M$). Hence it is also naturally identified 
with the set of conjugacy
classes of the elements of $\pi_1(CM)$. 
Since every contact manifold is
oriented and the contact structure was assumed to be cooriented, 
we get that the planes of
the contact structure are naturally oriented. This orientation induces the
orientation of the $S^1$-fibers of $\pr:CM\rightarrow M$. 
Put $f\in\pi_1(CM)$ to be the class of the oriented $S^1$-fiber of 
$\pr:CM\rightarrow M$.

The Theorem of Chow~\cite{Chow} and Rashevskii~\cite{Rashevskii} says
that every connected component of the space of curves contains a Legendrian 
curve. Straightforward verification shows that the insertion 
of the zig-zag into 
the Legendrian curve $K$ (see Figure~\ref{twocusp.fig}) 
changes the value of the $\spin$-structure on the corresponding framed curve. 
It is easy to verify (see~\cite{FuchsTabachnikov}) that the two connected components of the space of Legendrian 
curves that contain respectively $K$ and $K$ with the extra 
zigzag correspond to the 
conjugacy classes of $\vec K$ and of $\vec K f$ 
(or of $\vec K f^{-1}$) in $\pi_1(CM)$. (We 
obtain $\vec K f$ or $\vec K f^{-1}$ depending on which of the two possible
zig-zags we insert.)

Let $\mathcal L\subset \mathcal F$ be a connected component of the space of
Legendrian curves in $(M, C)$ that corresponds to the conjugacy class of
$\vec K \in \pi_1(CM)$. Then every connected component 
$\mathcal L'\subset \mathcal F$ of the space of Legendrian curves
corresponds to the conjugacy class of $\vec K f^{2n}\in\pi_1(CM)$,
for some $n\in\Z$.

Hence $\mathcal F$ satisfies condition~$(*)$ if and only if for every 
$n\neq 0$ the elements $\vec K$ and $\vec K f^{2n}$ are not conjugate in
$\pi_1(CM)$.

Assume that $\mathcal F$ does not satisfy condition~$(*)$, then there
exists a nonzero $n\in\Z$ and $\beta\in\pi_1(CM)$ such that 
\begin{equation}\label{conjugate}
\beta \vec K \beta^{-1}=\vec K f^{2n}\in\pi_1(CM, \vec K(1)).
\end{equation} 

This implies that $\pr_*(\beta)$ and $\pr_*(\vec K)$ commute in $\pi_1(M,
K(1))$. The commutation relation gives a mapping $\mu:T^2\rightarrow M$ of
the two-torus $T^2=S^1\times S^1$ such that $\mu |_{(1\times S^1)}=K$ and 
$\mu |_{(S^1\times 1)}=\pr(\beta)$. 

Put $e$ to be the Euler class of the contact bundle.
Consider the locally-trivial $S^1$-fibration $p:M'\rightarrow T^2$ induced by 
$\mu$ from the $S^1$-fibration $\pr :CM \rightarrow M$.
One can verify that $2n\in \Z=H^2(T^2, \Z)$ is the Euler class of
$p$. On the other hand the Euler class of $p$ is $\mu^*(e)$  and it is naturally 
identified with the value of $e$ on the homology class realized by $\mu(T^2)$. 
This implies that if $\mathcal F$ does not satisfy condition~$(*)$, then there 
exists a homology class $\alpha$ from the statement of the Proposition.

On the other hand the existence of the class $\alpha$ from the statement of
the Proposition implies that there exists a Legendrian curve $K\in\mathcal
F$ such that $\vec K$ is conjugate to $\vec K f^{n}$, for $n$ being the
value of $e$ (the Euler class of the contact bundle) on the homology class 
realized by $\mu(T^2)$. (Proposition~\ref{existcontact} says that 
$e=2\alpha$, for some $\alpha \in H^2(M, \Z)$, and 
hence $n$ is even.) This means that $\mathcal F$ does not satisfy
condition~$(*)$ and we have proved Proposition~\ref{interpretationconditionII}.
\qed

\subsection{Proof of Theorem~\ref{atoroidal}}\label{Proofatoridal}
\begin{emf}\label{proofconditionI}
Let $\pr:CM\rightarrow M$ be the locally trivial $S^1$-fibration introduced
in~\ref{h-principleLegendrian} and let $f\in\pi_1(CM)$ be the class of an
oriented $S^1$-fiber of $\pr$.

Similar to~\ref{ProofinterpretationconditionII} we get that to prove 
that all the components of the space of framed curves satisfy
condition~$(*)$, it suffices to show that
$\vec K$ and $\vec K f^{2n}$ are not conjugate in $\pi_1(CM)$, for all
$0\neq n\in\Z$ and $\vec K\in\pi_1(CM)$. 
Let $\vec K,\beta\in\pi_1(CM)$ and
$n\in\Z$ be such that 
\begin{equation}\label{conjugate1}
\beta \vec K \beta^{-1}=\vec K f^{2n}\in\pi_1(CM, \vec K(1)).
\end{equation} 
We have to show that $n=0$.

Proposition~\ref{commute} says that $f$ is in the center of $\pi_1(CM, \vec
K (1))$. Hence 
\begin{equation}\label{transpose}
\beta\vec K=\vec K\beta f^{2n}.
\end{equation}

Identity~\eqref{transpose} implies that $\pr_*(\beta)$ and $\pr_* (\vec K)$ 
commute in $\pi_1(M)$.
Hence there exists a mapping of the two-torus $\mu:T^2=S^1\times S^1\rightarrow M$ 
such that $\mu(S^1\times 1)=\pr (\vec K)$ and $\mu(1\times S^1)=\pr(\beta)$.
By the assumption of the Theorem
$\mu:\pi_1(T^2)\rightarrow\pi_1(M)$ has a nontrivial kernel. Thus there
exist $i,j\in\Z$ with at least one of $i$ and $j$ being nonzero such that 
$\pr(\vec K)^i=\pr(\beta)^j \in \pi_1(M)$, and hence
\begin{equation}\label{relatepower}
\vec K^i=\beta ^j f^l, \text{ for some } l\in \Z. 
\end{equation}

Since the situation is symmetric, we assume that $j\neq 0$.

Thus $\vec K^i \vec K^i= \vec K^i\beta ^j f^l=\beta ^j f^l\vec K^i$.
Applying~\eqref{transpose} to the last identity we get that $f^{2nij}=1$. 
Since $\pi_2(M)=0$ we see that $f$ has infinite order in $\pi_1(CM)$, and
hence $2nij=0$. If $n$ is zero, then we are done. Hence we have to look at
the case of $i=0$. (We assumed that $j\neq 0$.) From~\eqref{relatepower}
we get that $\beta ^j=f^{-l}$, and hence by Proposition~\ref{commute}
$\beta ^j$ is in the center of $\pi_1(CM)$. Thus $\beta^j \vec K=\vec K
\beta^j$. On the other hand using~\eqref{transpose} we get that 
$\beta ^j \vec K= \vec K \beta^j f^{2nj}$. Since $f$ has infinite order in
$\pi_1(CM)$ we get that $2nj=0$. By our assumptions $j\neq 0$ and we have
$n=0$.
This finishes the proof of Theorem~\ref{atoroidal}.\qed

\end{emf}
 
%
%
%
%
%
%

%
%
%
%
%
%
%
%

\subsection{Proof of Theorem~\ref{tight}}\label{Prooftight}
By Theorem~\ref{isomorphism} it suffices to show that all the connected components
of the space of framed curves in $M$ satisfy condition~$(*)$.

Let $e\in H^2(M,\Z)$ be the Euler class of the contact
bundle of $(M,C)$.
Proposition~\ref{interpretationconditionII} implies that 
it suffices to show that $e(\alpha)=0$, for every 
homology class $\alpha\in H_2(M)$ realizable by a mapping
$\mu:T^2\rightarrow M$.
The result of
D.~Gabai (see Corollary 6.18~\cite{Gabai}) implies that every 
$\alpha\in H^2(M, \Z)$ realizable by a mapping of $T^2$ 
can be realized by a collection of 
spheres and of a torus that are embedded into $M$. Finally the result of
Ya.~Eliashberg (see~\cite{Eliashbergtight} Theorem $2.2.1$) says that for a tight contact structure the
value of $e$ on any embedded torus or sphere is zero. Hence $e(\alpha)=0$. 
This finishes the Proof of Theorem~\ref{tight}.\qed

\subsection{Proof of Proposition~\ref{existcontact}}\label{Proofexistcontact}
{\em First we show that if $e\in H^2(M,\Z)$ can be realized as the Euler class of
the contact structure, then $e=2\alpha$ for some $\alpha\in H^2(M,\Z)$.\/} 

Since the contact structure is cooriented we get that the tangent bundle
$TM$ is isomorphic to the sum $C\oplus\epsilon$ of the oriented  contact bundle
$C$ and the trivial oriented line bundle $\epsilon$. 
The tangent bundle of every
orientable $3$-manifold is trivializable and we get that the second Stiefel-Whitney 
class of the contact bundle is zero. But the second Stiefel-Whitney
of $C$ is the projection of the Euler class of $C$ under the natural mapping
$H^2(M, \Z)\rightarrow H^2(M, \Z_2)$, and we get that $e=2\alpha$ for some
$\alpha\in H^2(M,\Z)$.

{\em Now we show that if $e=2\alpha$, for some $\alpha\in H^2(M,\Z)$, then $e$
can be realized as the Euler class of a cooriented contact structure on $M$.\/}

Consider an oriented $2$-dimensional vector
bundle $\xi$ over $M$ with the Euler class $e(\xi)=e=2\alpha\in H^2(M,\Z)$.
The second Stiefel-Whitney class $w_2(\xi)$ 
of $\xi$ is zero, since it is the projection of $e(\xi)=2\alpha\in H^2(M, \Z)$.
Since $\xi$ is an oriented vector bundle we have $w_1(\xi)=0$. 

Consider the sum $\xi\oplus\epsilon$ of $\xi$ with the trivial oriented
$1$-dimensional vector bundle $\epsilon$. Clearly the total Stiefel-Whitney class of 
the $3$-dimensional oriented vector bundle $\xi\oplus\epsilon$ is equal to $1$, and
the Euler class of $\xi\oplus\epsilon$ is equal to $0$. Using the
interpretation of the Stiefel-Whitney and the Euler classes of
$\xi\oplus\epsilon$
as obstructions for the trivialization of $\xi\oplus\epsilon$, we get that
$\xi\oplus\epsilon$ is trivializable. Since the tangent bundle of an oriented
$3$-dimensional manifold is trivializable, we see that $\xi$ is isomorphic 
to an oriented sub-bundle of $TM$. Since $M$ is oriented this sub-bundle of
$TM$ is also cooriented.
Now the Theorem of Lutz~\cite{Lutz}, that says that every homotopy 
class of distributions of $2$-planes tangent to $M$ contains a contact 
structure, implies the existence of a cooriented contact structure with the
Euler class $e$.\qed

\subsection{Proof of Theorem~\ref{example1}}\label{Proofexample1}
\begin{emf}\label{Proofaexample1}
{\em Proof of statement {\bf a} of Theorem~\ref{example1}.\/}
Clearly (see Figure~\ref{homotopy.fig}) the two Legendrian knots $K_0$ and
$K_1$ belong to the same component of the space of Legendrian curves. It is easy to see 
that if $K_0$ realizes the isotopy class of a framed knot $\tilde K^0$, then
$K_1$ realizes the isotopy class of $\tilde K^{-2}$ (see~\ref{obstruction} for
the definition of $\tilde K^{-2}$).
Below we show that $\tilde K^0$ and $\tilde K^{-2}$ are isotopic framed
knots.

\begin{figure}[htbp]
 \begin{center}
  \epsfxsize 10cm
  \hepsffile{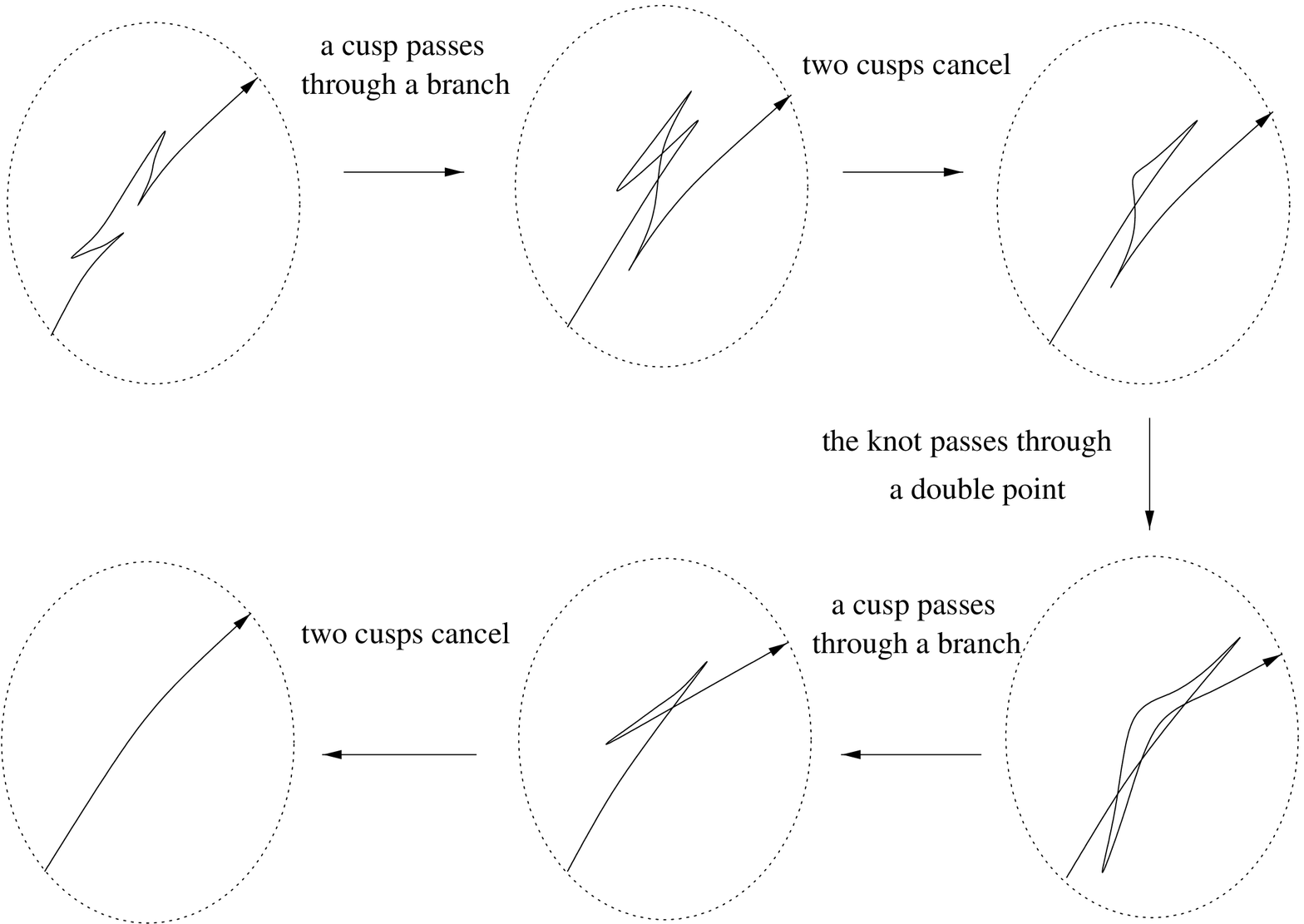}
 \end{center}
\caption{}\label{homotopy.fig}
\end{figure}

Let $t\times S^2\subset S^1\times S^2$ be the sphere that crosses $\tilde K^0$ at exactly one
point, and let $N=[0,1]\times S^2$ be a thin tubular neighborhood of $t\times
S^2$. Fix $x\in S^2$ (below called the North pole) and the direction in $T_x
S^2$ (below called the zero meridian). 
We can assume that the knot $\tilde K^0$ inside $N=[0,1]\times S^2$ looks as
follows: it intersects each $y\times S^2\subset N=[0,1]\times S^2$ at the
North pole of the corresponding sphere, and the framing of the knot is parallel to the zero meridian.

Consider an automorphism $\nu:S^1 \times  S^2\rightarrow S^1 \times  S^2$
that is identical outside of $N=[0,1]\times S^2$ such that it rotates each 
$y\times S^2\in [0,1]\times S^2$ by $4\pi y$ around the North pole in the clockwise 
direction. Clearly under this automorphism $\tilde K^0$ gets two extra negative 
twists of the framing and $\nu(\tilde K^0)=\tilde K^{-2}$. On the other hand 
it is easy to see that $\nu$ is diffeotopic to the identity, since it corresponds 
to the contractible loop in $SO(3)=\R P^3$. Hence we see that $\tilde K^0$ and 
$\tilde K^{-2}$ are isotopic framed knots. This finishes the proof of  
statement {\bf a} of Theorem~\ref{example1}.
\end{emf}

To prove statement {\bf b} of the Theorem we need the following
Proposition.

\begin{prop}\label{propS1xS2}
Let $C$ be a cooriented contact structure on $M=S^1\times S^2$
with a nonzero Euler class $e$ of the contact bundle. Let $CM$ be the spherical contact
bundle, let $\pr:CM\rightarrow M$ be the corresponding locally trivial 
$S^1$-fibration, and let $f\in\pi_1(CM)$ be the class of an oriented $S^1$-fiber 
of $\pr$. Then $f$ is of finite order in $\pi_1(CM)$ and $\pi_2(CM)=0$.
\end{prop}

\begin{emf}{\em Proof of Proposition~\ref{propS1xS2}.\/}
Consider the oriented $2$-plane bundle $p:\xi\rightarrow S^2$ that is the
restriction of the contact bundle over $M$ to the sphere 
$1\times S^2\subset S^1\times S^2$. The Euler class of $p$ is the value of 
$e$ on the homology class realized by $1\times S^2$, and hence is nonzero.
Let $S\xi$ be the manifold obtained by the fiberwise spherization of $p$, and
let $\bar p:S\xi\rightarrow S^2$ be the corresponding locally trivial
$S^1$-fibration. Since the Euler class of $p$ is nonzero we get that 
a certain multiple of the class of the fiber of $\bar p$ is homologous to
zero. But $\pi_1(S\xi)$ is generated by the class of the fiber, 
and hence the class of the fiber of $\bar p$ is of finite order in
$\pi_1(S\xi)$. This implies that $f\in\pi_1(CM)$ is of finite order.

The statement that $\pi_2(CM)=0$ follows from the exact homotopy sequence of 
$\pr:CM\rightarrow M$ and the fact that $f\in\pi_1(CM)$ is of finite order.
\qed
\end{emf}  

\begin{emf}{\em Proof of statement {\bf b} of Theorem~\ref{example1}.\/}
Let $\mathcal L$ be the connected component of the space of Legendrian
curves that contains $K_0$ and $K_1$.
Figure~\ref{homotopy.fig} shows that $K_0$ can be changed to $K_1$ (in the
space of Legendrian curves) by a
sequence of isotopies and one passage through a transverse double point. Hence if there 
exists a $\Z$-valued invariant $I$ of Legendrian knots from $\mathcal L$
that increases by one under every positive passage through a transverse double point of a Legendrian 
knot, then it distinguishes $K_0$ and $K_1$. (Clearly if such $I$ does
exist, then it is an order one invariant of Legendrian knots.) Below we show the existence of
such $I$ in the connected component $\mathcal L$.

Put $I(K_0)=0$. Let $K'\in\mathcal L$ be a Legendrian knot, and let 
$\gamma$ be a generic path in $\mathcal L$ connecting $K_0$ and $K'$. Let
$J_{\gamma}$ be the set of moments when $\gamma$ crosses the discriminant (i.e.
the subspace of singular knots) in $\mathcal L$, and let $\sigma_j$, $j\in
J_{\gamma}$, be the signs of these crossings. For a generic path $\gamma\subset \mathcal L$ put 
$\Delta_I(\gamma)=\sum_{j\in J_{\gamma}}\sigma_j$. It is clear that if $I$ (with
$I(K_0)=0$) does exist,
then $I(K')=\Delta_I(\gamma)$. To show that $I$ does exist we have to
verify that for every Legendrian knot $K'\in\mathcal L$ and for a generic
path $\gamma$ connecting $K'$ to $K_0$ the value of  
$\Delta_I(\gamma)$ does not depend on the choice of a generic
path $\gamma$ connecting $K_0$ and $K'$, or equivalently we have to show
that $\Delta_I(\gamma)=0$  for every generic closed loop $\gamma$ 
connecting  $K_0$ to itself. 

There are two codimension two strata of the discriminant of $\mathcal L$.
They are formed respectively by singular Legendrian knots with two 
transverse double points, and by Legendrian knots with one double point at
which the two intersecting branches are tangent of order one.
Straightforward verification shows that $\Delta_I(\beta)=0$, for every small closed loop
$\beta$ going around a codimension two stratum of $\mathcal L$. 

This implies that for every generic loop $\gamma$ connecting $K_0$ to itself 
the value of $\Delta_I(\gamma)$ depends only on the element of
$\pi_1(\mathcal L, K_0)$ realized by $\gamma$. Hence to prove the
existence of $I$ it suffices to show that $\Delta_I(\gamma)=0$ for every
$\gamma\in\pi_1(\mathcal L, K_0)$.

Clearly $\Delta_I(\gamma^p)=p\Delta_I(\gamma)$ and since $\Z$ is torsion
free, we get that to prove Theorem~\ref{example1} it suffices to show that 
for every $\gamma\in\pi_1(\mathcal L, K_0)$ there exists a nonzero $p\in \Z$
such that $\Delta_I(\gamma^p)=0$.

The $h$-principle says that the space of Legendrian curves in $(M,C)$ is weak homotopy equivalent to 
the space of free loops $\Omega CM$ in the spherical contact bundle $CM$ of 
$M=S^1\times S^2$. (The mapping giving the equivalence lifts a Legendrian
curve $K$ to a loop $\vec K$ in $CM$ by sending $t\in S^1$ to the point in
$CM$ that corresponds to the direction of the velocity vector of $K$ at $K(t)$.)

Thus $\pi_1(\mathcal L, K_0)$ is naturally isomorphic to $\pi_1(\Omega CM,
\vec K_0)$. Proposition~\ref{propS1xS2} says that 
$\pi_2(CM)=0$ and from~\ref{Hansen} we get that $\pi_1(\Omega
CM, \vec K_0)$ is isomorphic to the centralizer $Z(\vec K_0)$ of $\vec
K_0\in \pi_1(CM, \vec K_0(1))$. Using Propositions~\ref{commute}
and~\ref{propS1xS2} we see that either $\pi_1(CM)=\Z$ or
$\pi_1(CM)=\Z\oplus\Z_p$, for some nonzero $p\in\N$.
Hence there exists $n\in\Z$ and nonzero $m\in\Z$ 
such that $\gamma^m=\vec K_0^n\in\pi_1(CM, \vec K_0(1))$. (One should take
$n$ and $m$ to be divisible by $p$ if $\pi_1(CM)=\Z\oplus\Z_p$.)
But the loop
$\alpha$ in $\pi_1(\mathcal L, K_0)$ corresponding to $\vec K_0^n$ is just 
the sliding
of $K_0$ $n$ times along itself according to the orientation.
(This deformation is  induced by the rotation 
of the parameterizing circle.) This loop does not intersect the
discriminant, 
and hence $\Delta_I(\alpha)=0$. This finishes the proof of 
statement {\bf b} of Theorem~\ref{example1}.\qed
\end{emf}

\subsection{Proof of Theorem~\ref{example2}}\label{Proofexample2}
\begin{emf}\label{part1}{\em $K_1$ and $K_2$ are homotopic Legendrian
curves and they realize isotopic framed knots.\/}
Let $f_1\in\pi_1(CM)$ be the class of the $S^1$-fiber of the fibration 
$\pr:CM\rightarrow M$. The $h$-principle says that 
the connected component of the space of Legendrian curves that contains $K$
corresponds to the conjugacy class of $\vec K\in\pi_1(CM)$. From the formula
for the Maslov number deduced in~\cite{FuchsTabachnikov} 
it is easy to see that the connected components containing 
$K_1$ and $K_2$ correspond to the conjugacy classes of $\vec K_1=\vec K
f_1^r$ 
and  of $\vec K_2=\vec K f_1^{-r}$. 
Let $f_2\in\pi_1(CM)$ be an element projecting to the class $f\in\pi_1(M)$ of 
the $S^1$-fiber of $p:M\rightarrow F$. The value of the Euler class of the 
contact bundle  on the homology class realized by $\mu(T^2)$ is equal to 
$2r\in\Z$. (Here $\mu$ is the mapping from the description of the Euler
class of the contact bundle.) And because of the reasons explained in the proof of
Proposition~\ref{interpretationconditionII}
we get that $\vec K f_2=f_2 \vec K f_1^{2r}$, for the Legendrian knot $K$
used to construct $K_1$ and $K_2$. Now Proposition~\ref{commute}
implies that $\vec K_1$ and $\vec K_2$ are conjugate in $\pi_1(CM)$ and
hence $K_1$ and $K_2$ are in the same component of the space of Legendrian
curves. 

The fact that $K_1$ and $K_2$ realize isotopic framed knots is clear,
because as unframed knots they are the same, and as it is shown
in~\cite{FuchsTabachnikov} every pair of extra cusps corresponds to the
negative extra twist of the framing.
\end{emf}

\begin{emf}\label{ideaexample2}{\em The idea of the proof of the fact that $K_1$ and
$K_2$ can be distinguished by an order one invariant of Legendrian knots.\/}  
Let $d$ be a point in $M$. Let $K_s$ be a singular unframed knot with one 
double point.  The double point separates $K_s$ into two oriented loops. 
Deform $K_s$ preserving the double point, so that the double point is located
at $d$. Choosing one of the two loops of $K_s$ we obtain
an ordered set of two elements $\delta_1, \delta_2\in \pi_1(M,d)$, or which
is the same an element $\delta_1\oplus\delta_2\in\pi_1(M,
d)\oplus\pi_1(M,d)$. 
Clearly there is a unique element of the set $B$ that corresponds to the
original singular unframed knot $K_s$, where $B$ is is the quotient set of 
$\pi_1(M,d)\oplus\pi_1(M,d)$ modulo the consequent actions of the following groups:
\begin{description}
\item[1] $\pi_1(M)$ whose element $\xi$ acts on
$\delta_1\oplus\delta_2\in\pi_1(M)\oplus\pi_1(M)$ by sending it to 
$\xi\delta_1 \xi^{-1}\oplus \xi \delta_2
\xi^{-1}\in\pi_1(M)\oplus\pi_1(M)$.
(This corresponds to the ambiguity in deforming $K_s$, so
that the double point is located at $d$.) 

\item[2] $\Z_2$ that acts via the cyclic permutation of the two summands.
(This corresponds to the ambiguity in the choice of one of the two
loops of $K_s$.)
\end{description}

Thus we have a mapping $\nu$ from the set of singular unframed knots with 
one double point to $B$. Let $\alpha:B\rightarrow \Z$ be the function such
that 
\begin{description}
\item[a] $\alpha(b)=0$, provided that $b$ contains the class of
$1\oplus\delta\in\pi_1(M)\oplus\pi_1(M)$ for some $\delta\in\pi_1(M)$,
\item[b] $\alpha(b)=1$ otherwise.
\end{description}

Assume that $I^{\mathcal L}$ is an invariant of Legendrian knots 
from $\mathcal L$ such that under every (generic transverse) 
positive passage through a
discriminant in $\mathcal L$ it increases by $\alpha\circ\nu (K_s)$, where $K_s$ is the
unframed singular knot corresponding to the crossing of the discriminant.
Clearly such 
$I^{\mathcal L}$ is an order one invariant of framed knots from 
$\mathcal L$. To prove the Theorem we show the existence of such
$I^{\mathcal L}$, and then we show that it distinguishes 
$K_1$ and $K_2$.
\end{emf}

\begin{emf}\label{existence}{\em The existence of $I^{\mathcal L}$.\/}
Let $\gamma$ be a generic path in $\mathcal L$ that starts with $K_1$.
Let $J_{\gamma}$ be the set of instances when $\gamma$ crosses the discriminant (i.e.
the subspace of singular knots) in $\mathcal L$, and let $\sigma_j$, $j\in
J_{\gamma}$, be the
signs of these crossings. Let $J'_{\gamma}\subset J_{\gamma}$ be those instances for 
which the value of $\alpha\circ \nu$ on the corresponding singular unframed knots is $1$.
For a generic path $\gamma\subset \mathcal L$ put 
$\Delta^{\mathcal L}_I(\gamma)=\sum_{j'\in J'_{\gamma}}\sigma_{j'}$.

Similarly to~\ref{Proofexample1} we get that to prove the existence of $I^{\mathcal L}$
it suffices to show that $\Delta_I^{\mathcal L}(\gamma)=0$, for every generic closed loop $\gamma$.

Let $\mathcal C$ be the connected component of the space of unframed curves 
obtained by forgetting the framings on curves from $\mathcal F$, and let
$K'_1$ be the unframed knot obtained by forgetting the framing on $K_1\in
\mathcal L\subset \mathcal F$. Similarly to the above for a generic path 
$\gamma$ in $\mathcal C$ starting with $K'_1$ we put 
$\Delta^{\mathcal C}_I(\gamma)=\sum_{j'\in J'_{\gamma}}\sigma_{j'}$. (As above 
$J'_{\gamma}$ is the set of instances when the value of $\alpha\circ \nu$ on
the singular unframed knots obtained under $\gamma$ is equal to $1$, and $\sigma_{j'}$,
$j'\in J'_{\gamma}$, are the signs of the corresponding crossings of the discriminant.)

The codimension two stratum of the discriminant of $\mathcal C$ consists of 
singular curves whose only singularities are  two distinct transverse double 
points. Straightforward verification shows that $\Delta_I^{\mathcal C}(\beta)=0$
for any small loop $\beta$ going around the codimension two stratum. Hence
$\Delta_I^{\mathcal C}:\pi_1(\mathcal C,K'_1)\rightarrow \Z$ is a
homomorphism.

There are two codimension two strata in $\mathcal L$. They consist of
respectively singular Legendrian curves whose only singularities 
are two distinct transverse double points and of singular Legendrian curves
whose only singularity is one double point at which the two branches are
tangent. Considerations similar to the ones above show that 
$\Delta_I^{\mathcal L}:\pi_1(\mathcal L, K_1)\rightarrow \Z$ is a
homomorphism.

It is clear that if 
$\gamma'\in\pi_1(\mathcal C)$ is the element corresponding to 
$\gamma\in\pi_1(\mathcal L)$, then we have 
\begin{equation}\label{gammagamma'}\Delta_I^{\mathcal L}(\gamma)=
\Delta_I^{\mathcal C}(\gamma').
\end{equation}

The $h$-principle says that the space of Legendrian curves in $(M,C)$ is
weak homotopy equivalent to the space of free loops in the spherical
contact bundle  $CM$ of $M$. (The mapping that gives the equivalence lifts a
Legendrian curve $K$ in $(M, C)$ to a loop $\vec K$ in $CM$ by mapping $t\in
S^1$ to the point of $CM$ that corresponds to the velocity vector of $K$ at
$K(t)$.) Since $\pi_2(CM)=0$ for $M$ from the statement of the Theorem, 
we obtain (see~\ref{Hansen}) the natural isomorphism 
$t: \pi_1(\mathcal L, K_1)\rightarrow Z(\vec K_1)< \pi_1(CM, \vec K_1(1))$.
Since $\Delta^{\mathcal L}_I(\gamma^p)=p\Delta^{\mathcal L}_I(\gamma)$ and $\Z$ is torsion
free, we get that to show the existence of $I^{\mathcal L}$ it suffices to show 
that for every $\beta\in Z(\vec K_1)<\pi_1(CM,\vec K_1(1))$ there exist $0\neq n\in\Z$ 
and $\gamma\in\pi_1(\mathcal L, K_1)$ such that $t(\gamma)=\beta^n\in\pi_1(CM, \vec K_1(1))$ and
$\Delta^{\mathcal L}_I(\gamma)=0$.

Let $f\in\pi_1(M, K_1(1))$ be the class of the $S^1$-fiber of $p:M\rightarrow
F$. Let $f_1\in\pi_1(CM)$ be the class of an oriented $S^1$-fiber of
$pr:CM\rightarrow M$, and let $f_2$ be an element of $\pi_1(CM)$ such that 
$pr_*(f_2)=f\in\pi_1(M, K_1(1))$. 

Take $\beta\in Z(\vec K_1)$, then $\pr_*(\beta)\in Z(K_1)$.
Proposition~\ref{toughandtechnical} implies that there exist $0\neq n\in\Z$
and $i,j\in\Z$ such that $K_1^if^j=(\pr_*(\beta))^n\in\pi_1(M, K_1(1))$.
Using Proposition~\ref{commute} we get that 
\begin{equation}\label{expression}
\beta^n=\vec K_1^if_2^j f_1^l\text{ for some }i,j,l\in\Z. 
\end{equation}

As it was explained in~\ref{part1} we have 
\begin{equation}\label{extraobstruction}
\vec K_1 f_2= f_2 \vec K_1 f_1^{2r}.
\end{equation} 

Since $\beta\in Z(\vec K_1)$ we get that $\vec K_1\beta^n=
\beta^n \vec K_1$, and using~\eqref{expression} we see that 
$\vec K_1 \vec K_1^if_2^j f_1^l=\vec K_1^if_2^j f_1^l\vec K_1$.
Using~\eqref{extraobstruction}, Proposition~\ref{commute}, 
and the fact that $f_1$ has infinite order in $\pi_1(CM)$, we see that $j=0$
in~\eqref{expression}.

Hence 
\begin{equation}\label{newexpression}
\beta^n=\vec K_1^i f_1^l,\text{ for some }i,l\in\Z. 
\end{equation}

Clearly $f_1\in Z(\vec K_1)$ and since $t:\pi_1(\mathcal L, K_1)\rightarrow
Z(\vec K_1)$ is surjective
there exists a loop
$\gamma_3\in\pi_1(\mathcal L, K_1)$ such that $t(\gamma_3)=f_1$. 
Let $\gamma_2\in\pi_1(\mathcal L, K_1)$ be the loop corresponding to the
deformation under which $K_1$ slides once around itself according to the
orientation of $K_1$. (This deformation is induced by the rotation of the
circle parameterizing $K_1$.) Clearly $\gamma_2$ does not cross the
discriminant and hence $\Delta^{\mathcal L}_{I}(\gamma_2)=0$.

To prove the existence of $I^{\mathcal L}$ 
it suffices to show that $\Delta^{\mathcal L}_{I}(\gamma)=0$, for
$\gamma\in\pi_1(\mathcal L, K_1)$ such that $t(\gamma)=\vec K_1^i f_1^l$.
But this $\gamma$ is $\gamma_2^i \gamma_3^l$. Thus it
suffices to show that $0=\Delta^{\mathcal L}_{I}(\gamma)=
i\Delta^{\mathcal L}_{I}(\gamma_2)+l\Delta^{\mathcal L}_{I}(\gamma_3)$. Since 
$\Delta^{\mathcal L}_{I}(\gamma_2)=0$
we get that $\Delta^{\mathcal L}_{I}(\gamma)=l\Delta^{\mathcal L}_{I}(\gamma_3)$ and thus it suffices
to show that $\Delta^{\mathcal L}_{I}(\gamma_3)=0$.

Let $\gamma'_3\in\pi_1(\mathcal C, K'_1)$ be the loop corresponding to 
 $\gamma_3\in\pi_1(\mathcal L, K_1)$. Identity~\eqref{gammagamma'} says that 
$\Delta^{\mathcal L}_{I}(\gamma_3)=\Delta^{\mathcal C}_{I}(\gamma'_3)$.
Hence we have to show that $\Delta^{\mathcal C}_{I}(\gamma'_3)=0$.

In~\ref{explain} we show that $\gamma'_3\in\pi_1(\mathcal C, K'_1)$ can be realized 
as a power of the loop
described by the deformation shown in Figure~\ref{obstruction.fig}. 
Then since one of the loops of the only
singular knot arising under this deformation 
is contractible and the value of $\alpha\circ\nu$ on
such a singular knot is zero we get that $\Delta^{\mathcal
C}_{I}(\gamma'_3)=0$. This finishes the Proof of the existence of
$I^{\mathcal L}$ modulo the explanations given in~\ref{explain}. 
\end{emf}

\begin{emf}\label{explain}
{\em Now we show that $\gamma'_3\in\pi_1(\mathcal C, K'_1)$ can be realized 
as a sequence of loops
described by the deformation shown in Figure~\ref{obstruction.fig}.\/}
The $h$-principle for curves~\ref{h-principleforcurves} says that $\mathcal
C$ is weak homotopy equivalent to the space of free
loops $\Omega STM$ in the spherical tangent bundle $STM$ of $M$. In particular
$\pi_1(\mathcal C, K'_1)=\pi_1(\Omega STM, \vec K'_1)$.
In Subsubsection~\ref{Hansen} we introduced a
surjective homomorphism $t$ from $\pi_1(\Omega STM, \vec K'_1)$ onto
$Z(\vec K'_1)<\pi_1(STM, \vec K'_1(1))$. 

Let $\alpha, \beta\in\pi_1(\mathcal C, K'_1)$ be loops such that
$t(\alpha)=t(\beta)$. As it was explained
in~\ref{Hansen} the obstruction for $\alpha$ and $\beta$ to be homotopic
is an element of $\pi_2(STM)$. Since every orientable $3$-manifold is
parallelizable
we get that $STM=S^2\times M$. Clearly $\pi_2(M)=0$ for $M$ from 
the statement of the Theorem and
hence $\pi_2(STM)=\pi_2(S^2)=\Z$.

Consider the loop $\alpha'$ that looks the same as $\alpha$ except for a
small period of time when we perform the deformation shown in
Figure~\ref{obstruction.fig}. Clearly
$t(\alpha')=t(\alpha)=t(\beta)\in\pi_1(STM)$, and straightforward
verification show that the $\Z$-valued obstruction for $\alpha'$ and $\beta$
to be homotopic differs by one from the obstruction for $\alpha$ and
$\beta$ to be homotopic. Hence performing this operation (or its inverse)
sufficiently many times we can change $\alpha$ 
to be homotopic to $\beta$. 

Since $t(\gamma'_3)=t(1)=1\in\pi_1(STM)$ we get that
$\gamma'_3\in\pi_1(\mathcal C, K'_1)$ can be realized
as a power of the deformation shown in Figure~\ref{obstruction.fig}.

\begin{figure}[htbp]
 \begin{center}
  \epsfxsize 10cm
  \hepsffile{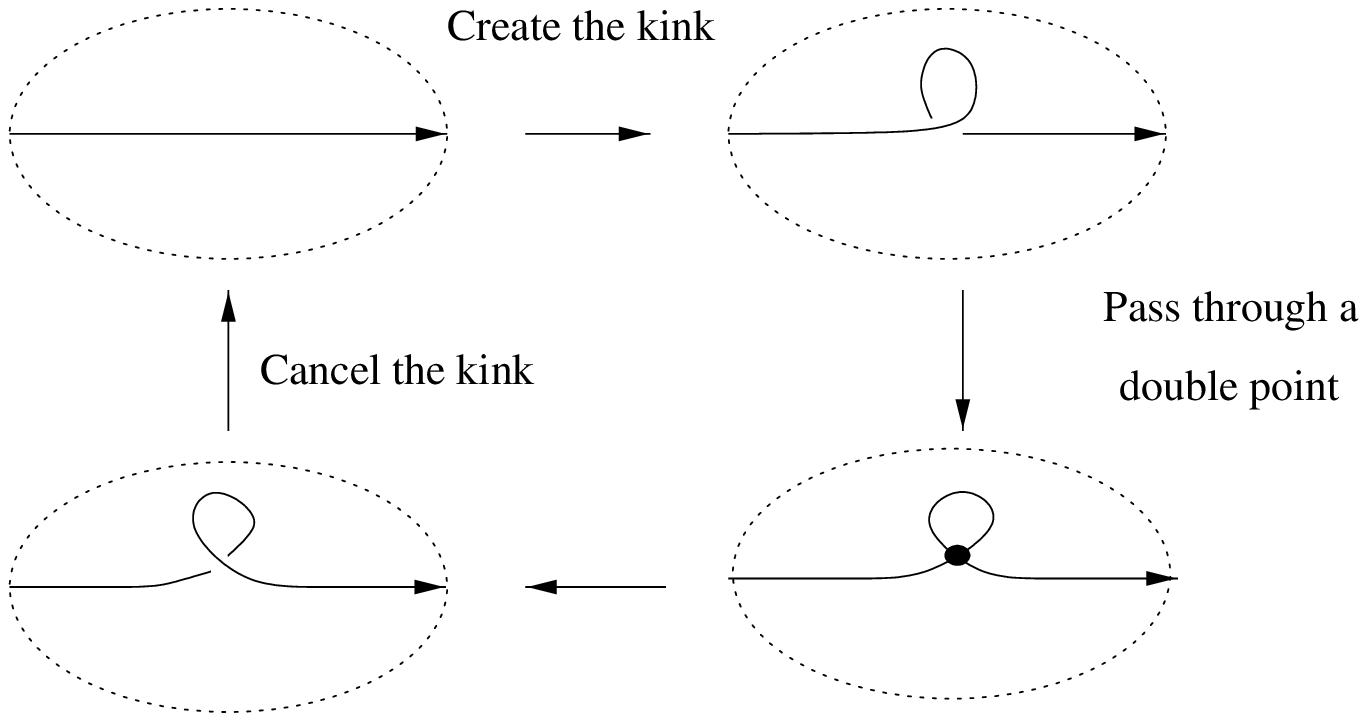}
 \end{center}
\caption{}\label{obstruction.fig}
\end{figure}

\end{emf}

\begin{emf}\label{distinguish}
{\em Let us show that $I^{\mathcal L}$ distinguishes $K_1$ and $K_2$.\/}
Let $\rho:[0,1]\rightarrow \mathcal L$ be a generic path connecting 
$K_1$ and $K_2$. To prove the Theorem we have to show that 
$\Delta_{I^{\mathcal L}}(\rho)\neq 0$.

Let $K'_1$ (resp $K'_2$) be the unframed knot obtained by forgetting the framing
on $K_1$ (resp $K_2$). 
Let $\rho ':[0,1]\rightarrow \mathcal C$ be the isotopy that deforms $K'_1$ into
$K'_2$ in
the category of unframed curves
under which $K'_1$ all the time 
stays in a thin tubular neighborhood of $K'_2$. Consider a homotopy  
$\bar \rho=S^1\rightarrow \mathcal C$ that corresponds to a product of paths
$\rho'\rho$. ($\bar \rho$ connects $K'_1$ to itself.) Clearly  
$\Delta_{I^{\mathcal C}}(\bar \rho)=\Delta_{I^{\mathcal L}}(\rho)$. 

For a loop $\alpha:S^1\rightarrow \mathcal C$ (that connects $K'_1$ to
itself) put $\lambda^{\alpha}:T^2=S^1\times S^1\rightarrow M$ to be a
mapping such that for every $t\in
S^1$ the mappings
$\lambda^{\alpha}\big|_{(t\times S^1)}$ and $\alpha(t) :S^1\rightarrow M$
are the same.
The value
of the Euler class of the contact bundle on the homology class realized by
$\lambda^{\bar \rho} (T^2)$
is equal to $2r\neq 0$.

Using the usual arguments we get that to prove the Theorem it suffices to 
show that for every $\gamma\in\pi_1(\mathcal C, K'_1)$ there exists 
$n\neq 0$ such that either 
\begin{description}
\item[a] the value of the Euler class of the contact bundle on the homology
class realized by $\lambda^{\gamma^n}:(T^2)\rightarrow M$ is zero,
or
\item[b] $\Delta_{I^{\mathcal C}}(\gamma^n)\neq 0$.
\end{description}

Consider the following loops $\gamma_1$ and $\gamma_2$.

{\em Loop $\gamma_1$.\/} Since $p(K'_1)$ is an orientation preserving loop
and $M$ is orientable, we get that the $S^1$-fibration over $S^1$
(parameterizing the knots) induced from $p:M\rightarrow F$ 
by $p\circ K'_1:S^1\rightarrow F$ is
trivializable.
Hence we can coherently orient the fibers of this fibration. The orientation
of the $S^1$-fiber over $t\in S^1$ induces the orientation of the
$S^1$-fiber
of $p$ that contains $K'_1(t)$.
The loop $\gamma_1$ is the deformation of $K'_1$
under which every point of $K'_1$ slides once around the fiber of $p$ that
contains this point
(staying inside the fiber) in the
direction specified by the orientation of the fiber corresponding to this
point.

{\em The loop $\gamma_2$\/} is the sliding of $K'_1$ along itself according
to the orientation.
(This deformation is induced by the rotation of the circle that
parameterizes
$K'_1$.)

The $h$-principle for curves says that $\mathcal C$ is weak homotopy
equivalent to the space of free loops $\Omega STM$ in the spherical tangent
bundle $STM$ of $M$. (The mapping that gives the equivalence lifts a
curve $K$ in $M$ to a loop $\vec K$ in $STM$ by mapping $t\in
S^1$ to the point of $STM$ that corresponds to the velocity vector of $K$ at
$K(t)$.) Let $t:\pi_1(\mathcal C, K'_1)\rightarrow Z(\vec K'_1)<\pi_1(\Omega STM, \vec K'_1)$
be the surjective homomorphism described in~\ref{Hansen}, and let $\bar
f\in\pi_1(STM)$
be the element that projects to the class $f\in\pi_1(M)$ of the
$S^1$-fiber of $p:M\rightarrow F$.

Using Proposition~\ref{toughandtechnical}
one verifies that for every $\gamma\in\pi_1(\mathcal C, K'_1)$ there exist $0\neq n\in\Z$ such 
$t(\gamma^n)=\bar f^i (\vec K'_1)^j=t(\gamma_1)^i t(\gamma_2)^j$, for some
$i,j\in\Z$. Let $\gamma_4$ be the loop described in Figure~\ref{obstruction.fig}.
(It is easy to see that $\gamma_4$ is in the center 
of $\pi_1(\mathcal C, K'_1)$.)
Similar to~\ref{explain} we get that $\gamma^n=\gamma_1^i \gamma_2^j\gamma_4^k$ 
for some $k\in\Z$.

It is easy to see that if $i=0$ then the value of the Euler class on the
homology
class realized by $\lambda^{\gamma^n}:T^2\rightarrow M$ is zero, and hence 
{\bf a} holds for $\gamma$. 

On the other hand if $i\neq 0$ then as we show below in~\ref{bholdsforgamma}
$\Delta_{I^{\mathcal
C}} (\gamma^n)\neq 0$ and {\bf b} holds for $\gamma$. (This finishes the
proof of the Theorem modulo the explanation below.)
\end{emf}

\begin{emf}\label{bholdsforgamma}
{\em If $i\neq 0$ then $\Delta_{I^{\mathcal
C}} (\gamma^n)\neq 0$ and {\bf b} holds for $\gamma$.\/}
The loop $\gamma_1$ crosses the discriminant twice, both crossings
occur with the same sign and the values of $\alpha\circ\nu$ on the corresponding
singular knots are equal to one. These crossings occur in the fiber over the double
point of $p(K'_1)$. (Since the double point of $p(K'_1)$ separates it into
two orientation reversing loops, the two points of $K'_1$ contained in 
this fiber induce opposite orientations of it, and the two branches of $K'_1$ 
that intersect the fiber slide in the opposite directions under $\gamma_1$.)
Hence $\Delta_{I^{\mathcal C}} (\gamma_1)=\pm 2\neq 0$. (The sign
depends on the orientation of the $S^1$-fibers of 
$T^2\rightarrow S^1$ used to induce the orientations of the fibers
containing the points of $K'_1$.)

Clearly 
$\Delta_{I^{\mathcal C}} (\gamma_2)=0$. 

Since one of the two loops of the only 
singular knot appearing in $\gamma_4$ is contractible, we have that
the value of $\alpha\circ \nu$ on the singular knot is $0$ and thus
$\Delta_{I^{\mathcal C}} (\gamma_4)=0$.

Hence if $i\neq 0$, then 
$\Delta_{I^{\mathcal C}} (\gamma^n)=i\Delta_{I^{\mathcal C}}
(\gamma_1)+j\Delta_{I^{\mathcal C}}(\gamma_2)+k\Delta_{I^{\mathcal
C}}(\gamma_4)=i\Delta_{I^{\mathcal C}} (\gamma_1)=i(\pm 2)\neq 0$ 
and hence {\bf b} holds for $\gamma$.

This finishes
the proof of Theorem~\ref{example2}.
\qed
\end{emf}

{\bf Acknowledgments.} 
I am very grateful to Stefan Nemirovski, Serge Tabachnikov and 
Oleg~Viro for the valuable discussions and
suggestions. I am deeply thankful to H.~Geiges and A.~Stoimenow for the 
valuable suggestions, and to O.~Baues, M.~Bhupal, 
N.~A'~Campo, A.~Cattaneo, J.~Fr\"ohlich, J.~Latschev, A.~Shumakovich, and V.~Turaev
for many valuable discussions.

This paper was written during my stay at the Max-Planck-Institut f\"ur
Mathematik (MPIM), Bonn, and it is a continuation of the research conducted at 
the ETH Zurich~\cite{Chernovpreprint}.
I would like to thank the Directors and the staff
of the MPIM and the staff of the ETH 
for hospitality and for providing the excellent working conditions.

\end{document}